\documentclass[reqno]{amsart}
\usepackage{amssymb,enumerate}
\usepackage{pstcol}
\usepackage{pst-plot,pst-node,graphicx}

\newcommand{\DefaultColors}{d}
\newcommand{\SlideColors}  {s}
\newcommand{\GrayColors}   {g}

\newcommand{\PaperSize}  {p}
\newcommand{\SlideSize}  {s}

\let\picsize\PaperSize

\let\palette\DefaultColors

\if\palette\DefaultColors
    \definecolor{math-green}         {cmyk}{0.40, 0.00, 0.50, 0.30}
    \definecolor{my-lawn-green}      {cmyk}{0.50, 0.00, 0.98, 0.20}
    \definecolor{math-red}           {cmyk}{0.00, 0.70, 0.70, 0.10}
    \definecolor{my-red}             {cmyk}{0.00, 0.90, 0.90, 0.00}
    \definecolor{math-blue}          {cmyk}{0.70, 0.70, 0.00, 0.10}
    \definecolor{my-blue}            {cmyk}{1.00, 1.00, 0.00, 0.10}
    \definecolor{my-light-blue}      {cmyk}{0.25, 0.25, 0.00, 0.05}
    \definecolor{my-purple}          {cmyk}{0.32, 0.82, 0.00, 0.06}
    \definecolor{my-orange}          {cmyk}{0.00, 0.18, 0.50, 0.00}
    \definecolor{my-goldenrod-yellow}{cmyk}{0.00, 0.00, 0.24, 0.00}
\else\if\palette\SlideColors
    \definecolor{math-green}         {cmyk}{0.40, 0.00, 0.50, 0.15}
    \definecolor{my-lawn-green}      {cmyk}{0.25, 0.00, 0.49, 0.00}
    \definecolor{math-red}           {cmyk}{0.00, 0.50, 0.50, 0.00}
    \definecolor{my-red}             {cmyk}{0.00, 0.50, 0.50, 0.00}
    \definecolor{math-blue}          {cmyk}{0.50, 0.50, 0.00, 0.00}
    \definecolor{my-blue}            {cmyk}{0.5, 0.5, 0.00, 0.00}
    \definecolor{my-light-blue}      {cmyk}{0.12, 0.12, 0.00, 0.0}
    \definecolor{my-purple}          {cmyk}{0.16, 0.41, 0.00, 0.00}
    \definecolor{my-orange}          {cmyk}{0.00, 0.09, 0.25, 0.00}
    \definecolor{my-goldenrod-yellow}{cmyk}{0.00, 0.00, 0.12, 0.00}
\else\if\palette\GrayColors
    \definecolor{math-green}   {gray}{0.0}
    \definecolor{my-lawn-green}{gray}{0.7}
    \definecolor{math-red}     {gray}{0.0}
    \definecolor{my-red}       {gray}{0.0}
    \definecolor{math-blue}    {gray}{0.0}
    \definecolor{my-blue}      {gray}{0.5}
    \definecolor{my-light-blue}{gray}{0.85}
    \definecolor{my-purple}    {gray}{0.3}
    \definecolor{my-orange}    {gray}{0.8}
    \definecolor{my-goldenrod-yellow}{gray}{0.95}
\fi\fi\fi

\definecolor{light-gray}{gray}{0.8}

\newlength\braceln 

\makeatletter

\def\@@enum@[#1]{%
  \@enLab{}\let\@enThe\@enQmark
  \@enloop#1\@enum@
  \expandafter\edef\csname label\@enumctr\endcsname{\the\@enLab}%
  \expandafter\let\csname the\@enumctr\endcsname\@enThe
  \csname c@\@enumctr\endcsname7
  \expandafter\settowidth
            \csname leftmargin\romannumeral\@enumdepth\endcsname
            {\the\@enLab\hspace{\labelsep}}%
  \@enum@}

 \newcommand{\df}[1]{{\rmfamily\itshape\mdseries#1}}

 \let\paragraphname\@empty  

\def\subsection{\@startsection{subsection}{2}%
  \z@{.5\linespacing\@plus.7\linespacing}{-.5em}%
  {\normalfont\sffamily\bfseries\upshape}}

\renewenvironment{description}{\list{}{%
  \setlength\leftmargin{1.5em} \setlength\itemindent{-1.5em}
  \labelwidth\z@ }%
}{
  \endlist
}

\numberwithin{equation}{section} 

\newcommand{\customqed}[1]{{\renewcommand{\qedsymbol}{#1}\qed}}
\newcommand{\varqed}{\customqed{\hbox{$\diamondsuit$}}}

\newtheoremstyle{sfdefinition}
  {5pt}
  {5pt}
  {}
  {}
  {\sffamily\mdseries\upshape}
  {.}
  {.5em}
  {}

\newtheoremstyle{sftheorem}
  {5pt}
  {5pt}
  {\itshape}
  {}
  {\sffamily\mdseries\upshape}
  {.}
  {.5em}
  {}

 \theoremstyle{sftheorem}

 \theoremstyle{sftheorem}

 \newtheorem{lemma}{Lemma}[section]

 \newtheorem{theorem}[lemma]{Theorem}

 \theoremstyle{sfdefinition}

 \theoremstyle{remark}

 \newtheorem{Condition}[lemma]{Condition}

 \newenvironment{condition}{%
   \begin{Condition}}{\varqed\end{Condition}}

 \newtheorem{Example}[lemma]{Example}
 \newtheorem{Examples}[lemma]{Examples}
 
 \newenvironment{example}{%
   \begin{Example}}{\varqed\end{Example}}

 \newtheorem{Remark}[lemma]{Remark}
 \newtheorem{Remarks}[lemma]{Remarks}
 
 \newenvironment{remark}{%
   \begin{Remark}}{\varqed\end{Remark}}
 \newenvironment{remarks}{%
   \begin{Remarks}}{\varqed\end{Remarks}}

\newenvironment{alphenum} 
 {\begin{enumerate}[{\upshape (a)}]}{\end{enumerate}}

 {\begin{enumerate}[{\upshape (a)}]}{\end{enumerate}}

 \newcommand\ol[1]{\overline{#1}}
 \newcommand\ul[1]{\underline{#1}}
 
 \newcommand \ag{\alpha}
 \newcommand \bg{\beta}
 
 \newcommand \Gg{\Gamma}

 \newcommand \eps{\varepsilon}

 \renewcommand \lg{\lambda} 
 \newcommand \Lg{\Lambda}

 \newcommand \sg{\sigma}

 \newcommand{\cA}{\mathcal{A}}
 \newcommand{\cB}{\mathcal{B}}
 \newcommand{\cC}{\mathcal{C}}
 
 \newcommand{\cE}{\mathcal{E}}
 \newcommand{\cF}{\mathcal{F}}
 \newcommand{\cG}{\mathcal{G}}
 \newcommand{\cH}{\mathcal{H}}
 \newcommand{\cI}{\mathcal{I}}

 \newcommand{\cM}{\mathcal{M}}

 \newcommand{\cQ}{\mathcal{Q}}
 
 \newcommand{\cS}{\mathcal{S}}

 \newcommand{\cV}{\mathcal{V}}
 \newcommand{\cW}{\mathcal{W}}

 \newcommand{\MM}{\mathbb{M}}

 \newcommand{\ZZ}{\mathbb{Z}}

 \newcommand\setof[1]{\mathopen\{\,#1\,\mathclose\}}
 
 \newcommand\ang[1]{{\langle #1\rangle}}
 \newcommand\bigparen[1]{{\bigl(\,#1\,\bigr)}}

 \newcommand\flo[1]{{\lfloor #1\rfloor}}
 
 \newcommand\Flo[1]{{\left\lfloor #1\right\rfloor}}

 \newcommand\Prob{\mathop{\operator@font Prob}}

\newcommand{\half}{\frac{1}{2}}

  \renewcommand{\le}{\leqslant}
  \renewcommand{\ge}{\geqslant}
 
 \let\et=\wedge

 \newcommand{\sbs}{\subset}
 \newcommand{\sbsq}{\subseteq}
 \newcommand{\sps}{\supset}

 \newcommand{\xcpt}{\mathbin{\raise0.15ex\hbox{$\smallsetminus$}}}
 \newcommand{\eqdef}{\stackrel{\raise0.2ex\hbox{\scriptsize def}}{=}}

 \newcommand{\txt}[1]{\text{\rmfamily\mdseries\upshape{#1}}}

\makeatother

 \renewcommand{\d}{d} 
 \newcommand{\f}{f}
 \newcommand{\g}{g}
 \newcommand{\h}{h}
 \newcommand{\p}{p}
 
 \newcommand{\ncln}{\overline q}
 \renewcommand{\r}{r} 
 \newcommand{\s}{s}

 \newcommand{\aux}{c}
 \newcommand{\fxp}{\phi}
 \newcommand{\gxp}{\gamma}
 \newcommand{\txp}{\tau}
 \newcommand{\hxp}{\chi}

 \newcommand{\Bvalues}{\text{\rmfamily\mdseries\upshape Bvalues}}
 
 \newcommand{\Hvalues}{\text{\rmfamily\mdseries\upshape Hvalues}}
 \newcommand{\Btypes}{\text{\rmfamily\mdseries\upshape Btypes}}
 
 \newcommand{\Htypes}{\text{\rmfamily\mdseries\upshape Htypes}}

 \newcommand{\slope}{\text{\rmfamily\mdseries\slshape slope}}
 \newcommand{\minslope}{\text{\rmfamily\mdseries\slshape minslope}}
 \newcommand{\Type}{\text{\rmfamily\mdseries\slshape Type}}
 \newcommand{\Rank}{\text{\rmfamily\mdseries\slshape Rank}}
 \newcommand{\Body}{\text{\rmfamily\mdseries\slshape Body}}
 \newcommand{\sizelb}{\underline \Delta}
 \newcommand{\sizeub}{\overline \Delta}

 \newcommand{\bub}{{\Delta}}
 \newcommand{\bubxp}{\delta}

 \newcommand{\slb}{{\sigma}}

 \newcommand{\pub}{{\overline \p}}
 
 \newcommand{\R}{R}
 \newcommand{\Rmax}{\overline R}
 \newcommand{\T}{T}

 \newcommand{\dlb}{{\underline d}}
 \newcommand{\dub}{{\overline d}}
 \newcommand{\elb}{{\underline e}}
 \newcommand{\eub}{{\overline e}}

 \newcommand{\slopeincr}{\Lg}

\begin{document}

  \title[Compatible sequences]
       {Compatible sequences 
\\ and a slow Winkler percolation}

  \author{Peter G\'acs}\thanks{The paper was written during the author's
visit at CWI, Amsterdam, partially supported by a grant of NWO}
  \address{Boston University}
  \email{gacs@bu.edu}
  \date{\today}

 \begin{abstract}
Two infinite 0-1 sequences are called \df{compatible} when it is
possible to cast out 0's from both in such a way that 
they become complementary to each other.
Answering a question of Peter Winkler, we show that if the two 
0-1-sequences 
are random i.i.d. and independent from each other, with probability $\p$ of
1's, then if $\p$ is sufficiently small
they are compatible with positive probability.
The question is equivalent to a certain dependent percolation
with a power-law behavior: the probability that the origin is blocked at
distance $n$ but not closer decreases only polynomially
fast and not, as usual, exponentially.
 \end{abstract}

 \maketitle

 \section{Introduction}

 \subsection{The model}

Let us call any strictly increasing sequence $t = (t(0)=0, t(1),\dotsc)$ of
integers a \df{delay sequence}.
For an infinite sequence $x=(x(0), x(1),\dotsc)$, the 
delay sequence $t$ introduces a timing arrangement in which 
the value $x(n)$ occurs at time $t(n)$.
For two infinite 0-1-sequences $x_{\d}$
($\d = 0, 1$) and corresponding delay sequences $t_{\d}$ we say that there
is a \df{collision} at $(\d,n)$ if $x_{\d}(n) = 1$, and there is no $k$
such that $x_{1-\d}(k)=0$ and $t_{\d}(n)=t_{1-\d}(k)$.
We say that the sequences $x_{\d}$ are \df{compatible} if there is a pair
of delay sequences $t_{\d}$ without collisions.
It is easy to see that this is equivalent to saying that 0's can be
deleted from both sequences in such a way that the resulting sequences have
no collisions in the sense that they never have a 1 in the same position.

 \begin{example}
The sequences
 \[
  \begin{split}
      & \mathtt{0001100100001111\dots},
\\    & \mathtt{1101010001011001\dots}
  \end{split}
 \]
are not compatible.
The sequences $x,y$ below, are.
(We insert a $\mathtt{\ul{1}}$ in $x$, 
instead of deleting the corresponding $\mathtt{0}$ of $y$.)
 \[
\begin{split}
     x    &= \mathtt{0000100100001111001001001001001\dots},
\\   y    &= \mathtt{0101010001011000000010101101010\dots},
\\   x'   &= 
\mathtt{00001001\ul{1}000011110010\ul{1}0100100     1001\dotsb},
\\   y'   &= 
\mathtt{010101000     1011000000010     101101\ul{1}010\dotsb}.
\end{split}
\]
 \end{example}

Suppose that for $\d = 0, 1$, $X_{\d} = (X_{\d}(0),X_{\d}(1),\dotsc)$
are two independent infinite sequences of independent random variables
where $X_{\d}(j) = 1$ with probability $\p$ and 0 with probability $1-\p$.
Our question is: are $X_{0}$ and $X_{1}$ compatible with positive
probability?
The question depends, of course, on the value of $\p$: intuitively, it
seems that they are compatible if $\p$ is small, and our result will
confirm this intuition.

We can interpret the sequences $X_{\d}(j)$ as a chat of two members
of a retirement home with
unlimited time on their hands for maintaining a somewhat erratic
conversation.
Each member can be speaking or listening in any of the time periods 
$[j, j + 1)$ (value $X_{\d}(j) = 1$ or 0).
There is a nurse who can put a member to sleep for any of these time
periods or wake him up.
Speaker $\d$ wakes up to his $n$th action at time $t_{\d}(n)$.
The nurse wants to arrange that every time
when one of the parties is speaking the other one is up and listening.
She is really a fairy since she is clairvoyant: she sees the pair of
infinite sequences $X_{\d}$ (the talking/listening decisions) in advance
and can tailor her strategy $t_{\d}$ to it.
Here are some
interpretations that are either less frivolous or more hallowed by the
tradition of distributed computing.
 \begin{description}

  \item[Communication]
Assume that the sequences $X_{0}, X_{1}$ belong to two processors.
Value $X_{0}(i)=1$ means that in its $i$th turn, processor 0 wants to
send some message to processor 1, and $X_{0}(i)=0$ means that it is willing
to receive some message.
Assume that each processor is allowed to send an extra message 
(or, just to stay idle, to ``skip turns'')
in any time period, postponing the rest of its actions.
The goal is to achieve that with
the new sequences, whenever processor $\d$ is sending a message,
processor $1-\d$ is listening.

  \item[Dining]
Consider the following twist on the ``dining philosophers'' problem
(see~\cite{DijkstraDining71}).
Two philosophers, 0 and 1, sit across a table, with
a fork on both sides between them.
The philosophers have two possible actions: thinking and eating.
Sequence $X_{\d}(i)=1$ says that in her $i$th turn, philosopher $\d$ wants
to eat.
Both philosophers need two forks to eat, so any one will only be able to eat
when the other one is thinking.
Assume that both philosophers can be persuaded to insert some extra
eating periods into their sequences (or, as an equivalent but less
decorous possibility, to delete some thinking periods).

  \item[Queues]
A single server serves two queues, numbered by 0 and 1, where the queue
$\d$ of requests is being sent by a single user $\d$.
In both queues, a sequence of requests is coming in at discrete times
0,1,2, $\dotsc$.
Let $\tau_{\d}(i)$ be the time elapsed between the $i$th and $(i+1)$th
request in queue $\d$.
All variables $\tau_{\d}(i)$ are independent of each other, with
$\Prob\setof{\tau_{\d}(i)=n} = \p(1-\p)^{n-1}$ for $n > 0$
(discrete approximation of a Poisson arrival process).
Suppose that the senders of the queues will accept
\emph{faster service}: they are willing to send, instead of the originally
planned sequences $\tau_{\d}$, new
sequences $\tau'_{0}, \tau'_{1}$  where $\tau'_{\d}(i) \le \tau_{\d}(i)$ 
for all $\d,i$.
We want these new sequences to be served 
simultaneously by the single server.
Equivalently, suppose that both senders are willing to accept extra
services inserted into their queues.

 \end{description}

The question in all three cases is whether a scheduler who knows both
infinite sequences in advance, can make the needed synchronizations with
positive probability.

Peter Winkler and Harry Kesten, independently of each other, found an upper
bound smaller than $\half$ on the values $\p$ for which $X_{0}, X_{1}$ are
compatible.
We reproduce here informally Winkler's argument; 
it would be routine to formalize it.
Suppose that 
the infinite sequences $X$ and $Y$ are compatible, and let us denote by 
$X^{k},Y^{k}$ their $k$th initial segments.
Then we can delete some 0's from these segments in such a way that one of 
the resulting finite sequences, $X',Y'$, is
the complement of a prefix of the other: say, $X'$ is the complement of
a prefix of $Y'$.
Assume that both sequences contain at least $k/2 - \eps k$ 1's.
Then the number of deleted 0's in both sequences can be at most $2\eps k$.
Then we can reproduce the pair of sequences $X^{k},Y^{k}$ using the
following information.
 \begin{enumerate}[\upshape 1.]
   \item The sequence $Y'$.
   \item A 0-1 sequence $u$ of length $k$ whose 1's show the
positions of the deleted 0's in $X^{k}$.
   \item A 0-1 sequence $v$ of length $k$ whose 1's show the
positions of the deleted 0's in $Y^{k}$.
 \end{enumerate}
Using $Y'$ and $u$, we can restore $X^{k}$; using $Y'$ and $v$,
we can restore $Y^{k}$.
For $h(\p) = -\p\log_{2}\p - (1-\p)\log_{2}(1-\p)$, the
total entropy of the three sequences is at most $k + 2 k h(2\eps)$.
But the two sequences $X^{k},Y^{k}$ which we restored have total entropy
$2k$.
This gives an implicit lower bound on $\eps$: $h(2 \eps) \ge \half$.

Computer simulations by John Tromp
suggest that when $\p < 0.3$, with positive probability the sequences are
compatible.
The following theorem answers a question of Peter Winkler.

  \begin{theorem}[Main]\label{t.main}
If $\p$ is sufficiently small then with positive probability, $X_{0}$ and
$X_{1}$ are compatible.
  \end{theorem}

The threshold for $\p$ obtained from the proof is only $10^{-246}$, so
there is lots of room for improvement between this number and the
experimental $0.3$.
It turns out that, when deciding whether to cast out a 0 in position $n$ of
a sequence, we do not have to look ahead further than position $2n^{1.5}$
(see Subsection~\ref{ss.lookahead}).

  \subsection{A percolation}\label{ss.perc}
Compatibility can also be defined for finite sequences, and we can ask then
whether there is a polynomial 
algorithm that, given sequences $X_{0},X_{1}$ 
of length $n$, decides whether they are compatible.
It is easy to recognize that a dynamic programming algorithm will do it,
and that the structure created by the dynamic programming leads to a useful
reformulation of the original problem, too.

We define a directed graph $G = (V,E)$ as follows.
$V=\ZZ_{+}^{2}$ is the set of points $(i,j)$ where $i,j$ are nonnegative
integers.
When representing the set $V$ of points $(i,j)$ graphically, the right
direction is the one of growing $i$, and the upward direction is the one of
growing $j$.
The set $E$ of edges consists of all pairs of the form $((i,j),(i+1,j))$,
$((i,j),(i,j+1))$ and $((i,j),(i+1,j+1))$.

\begin{figure}
\if\picsize\SlideSize
    \psset{unit=0.5cm}
\else\if\picsize\PaperSize
    \psset{unit=0.7cm}
\fi\fi
        
\if\palette\DefaultColors
    \definecolor{my-blue-violet}{cmyk}{0.60, 0.60, 0.00, 0.10}
\else\if\palette\SlideColors
    \definecolor{my-blue-violet}{cmyk}{0.30, 0.30, 0.00, 0.00}
\else\if\palette\GrayColors
    \definecolor{my-blue-violet}{gray}{0.3}
\fi\fi\fi

\includegraphics{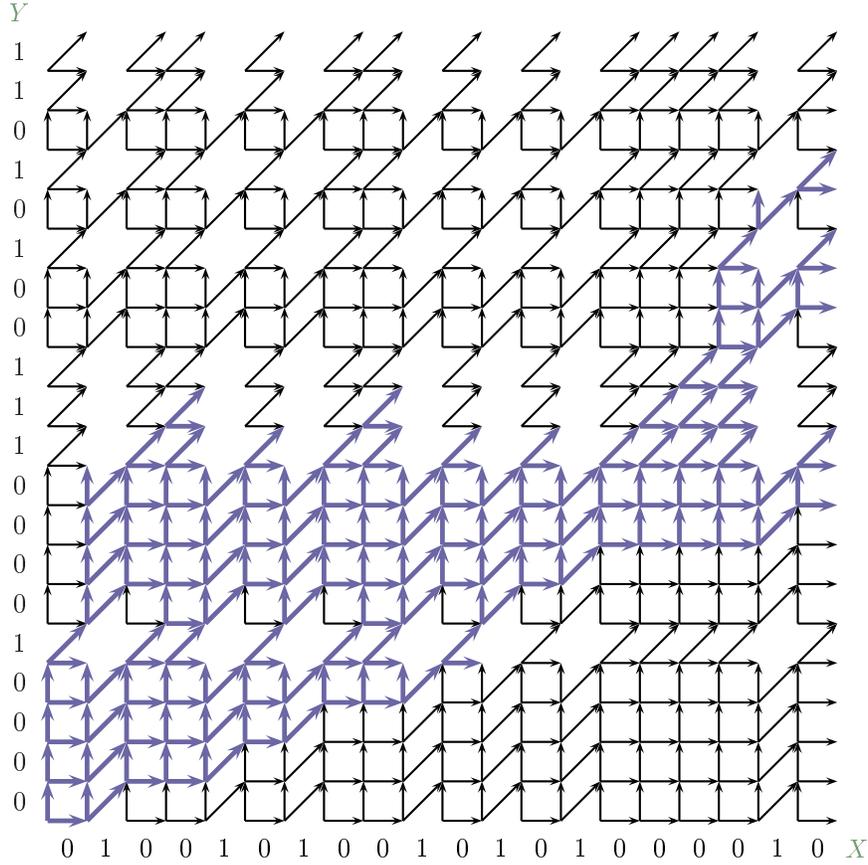}
 \caption{ \label{fig.chat-graph}
Oriented percolation arising from the compatibility problem}
\end{figure}

Sometimes we will write $X(i)$ for $X_{0}(i)$ and $Y(i)$ for
$X_{1}(i)$.
In the chat interpretation, when $X(i) = 1$ then
participant $0$ wants to speak in the $i$th turn of his
waking time, which is identified with the interval $[i, i+1)$.
In this case, we erase all edges of the form $((i,j),(i+1,j))$ for all
$j$ (this does not allow participant $1$ to sleep through this interval).
Similarly, when $Y(i) = 1$ then participant $1$ wants to speak in the
$i$th turn, and we erase all edges of the form $((j, i), (j, i+1))$.
If $X(i) = Y(j)$ then we also erase edge $((i, j), (i+1, j+1))$.
For $X(i) = 1$ since
we do not allow the two participants to speak simultaneously,
and for $X(i) = 0$ since the edge is not needed anyway, and this will allow
a nicer mathematical description.
This defines a graph $G(X, Y)$.
For an example, see Figure~\ref{fig.chat-graph}.
It is now easy to see that $X$ and $Y$ are compatible if and only if the
graph $G(X, Y)$ contains an infinite path starting at $(0,0)$.
We will say that \df{there is percolation} for $p$ if the probability at
the given parameter $\p$ that there is an infinite path is positive.
We will also use graph $G(X,Y)$ for the case of finite sequences
$X(0),\dots,X(m-1)$, $Y(0),\dots,Y(n-1)$, over $[0, m] \times [0, n]$.

We propose to call this sort of percolation, where two infinite random
sequences $X,Y$ are given on the two coordinate axes and the openness of a
point or edge at position $(i,j)$ depends on the pair $(X(i),Y(j))$, a
\df{Winkler percolation}.
Other examples will be seen in Section~\ref{s.related}.
Since we will talk about reachability a lot, the following notation is
useful.
If $u,v$ are points of a directed
graph (and the graph itself is clearly given from
the context), then 
 \[
   u \leadsto v
 \]
denotes the fact that $v$ is reachable from $u$ on a directed path.

 \subsection{Power-law behavior}\label{ss.power-law}

Our percolation problem has a
power-law behavior: the probability that the origin is blocked at
distance $n$ but not closer decreases only polynomially
fast, not, as usual, exponentially.

\begin{figure}

  \includegraphics{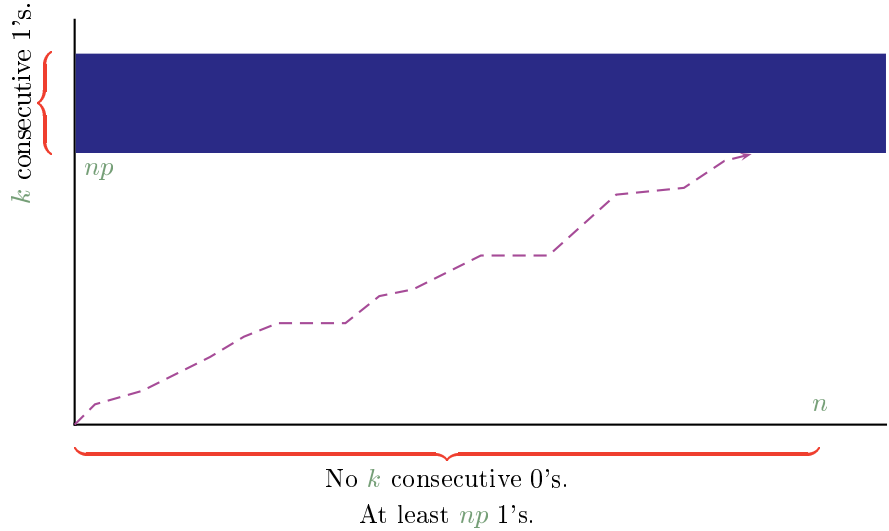}
  
 \caption{ \label{fig.power-law}
Power-law behavior.
For $k=c_{1}\log n$, the probability that
the above three events hold simultaneously is
$\ge n^{c_{2}(\log p)/p}$.}
 \end{figure}

Let us indicate informally the reason: a formal proof for
a similar problem is given in \cite{GacsClairv99}.
Let $W(n,k)$ be the event that $X(i) = 1$ for $i$ in $[n, n+k-1]$.
We can view this event as the occurrence of a vertical ``wall'' of width
$k$ at position $n$.
Let $\mu_{k}$ be the first number $n$ with $W(n, k)$.
Let $H(n, k)$ be the event that $Y(i) = 0$ for $i$ in $[n, n+k-1]$.
We can view this event as the occurrence of a horizontal ``hole'' of width
$k$ at position $n$.
Let $\nu_{k}$ be the first number $n$ with $H(n, k)$.
For a given $k$, $n = e^{k \p}$, let
 \[
  \p_{k} = \Prob\setof{\mu_{k} < \p n < n < \nu_{k}}.
 \]
Then with probability $\p_{k}$, a vertical wall occurs at position $\p n$
but no horizontal hole appears up to height $n$.
We can also assume that $\sum_{i \le n} Y(i) > \p n$.
Therefore every path will have to move right $> \p n$ steps while
ascending to height $n$, but then it will hit the wall which no hole can
penetrate up to height $n$.
This gives blocking at distance $n$ with approximately probability $p_{k}$.
See Figure~\ref{fig.power-law}.
 
It can be computed that 
this probability is a polynomial function of $n$ with a degree independent
of $k$.

\section{Outline of the proof}\label{s.outline}

\subsection{Renormalization}\label{ss.renorm}

The proof method used is \df{renormalization}, or
\df{multi-scale analysis}, and it is used frequently in statistical
mechanics.
The method is messy, laborious, and rather crude (rarely suited
to the computation of exact constants).
However, it is robust and well-suited to ``error-correction'' situations.
Here is a rough first outline.

 \begin{enumerate}[\upshape 1.]

  \item
Fix an appropriate sequence 
$\bub_{1} < \bub_{2} < \dotsb$, of scale parameters 
with $\bub_{k+1} > 4\bub_{k}$. 
Let 
 \[
   \cF_{k}
 \]
be the event that point $(0,0)$ is blocked in $[0, \bub_{k}]^{2}$.
(In other applications, it could be some other \df{ultimate bad event}.)
We want to prove 
 \[
   \Prob(\bigcup_{k} \cF_{k}) < 1.
 \]
This will be sufficient: if $(0,0)$ is not blocked in any
finite square then by compactness (or by what is sometimes called K\"onig's
Lemma), there is an infinite path starting at $(0,0)$.

  \item
Identify some events that you may call
\df{bad events} and some others called \df{very bad events},
where the latter are much less probable.
 
  \item\label{i.Fprime}
Define a series $\cM^{1},\cM^{2},\dotsc$ of models
similar to each other, where the very bad events of $\cM^{k}$ become
the bad events of $\cM^{k+1}$.
Let 
 \[
   \cF'_{k}
 \]
hold iff some bad event of $\cM^{k}$ happens in $[0, \bub_{k+1}]^{2}$.
  \item
Prove
 \begin{equation}\label{e.bd-non-reach-by-wall-all}
  \cF_{k} \sbs \bigcup_{i \le k} \cF'_{i}.
 \end{equation}

  \item
Prove $\sum_{k} \Prob(\cF'_{k}) < 1$.
  
  \end{enumerate}
In later discussions, we will
frequently delete the index $k$ from $\cM^{k}$ as well as from other
quantities defined for $\cM^{k}$.
In this context, we will refer to $\cM^{k+1}$ as $\cM^{*}$.

\subsection{Application to our case}

In our setting, the role of the
``bad events'' of Subsection~\ref{ss.renorm} will be played by \df{walls}.
In Subsection~\ref{ss.power-law}, we have seen a simple kind of wall coming
from a block of 1's.
Another kind of wall comes from two blocks of 1's that are too close to
each other.
The same example extrapolates to a more complicated kind of wall: let
$E(n, k, m)$ be the event that no hole $H(i, k)$ occurs on the segment
$[n, n+m]$.
For really large $m$, this event creates a horizontal wall.
Indeed, if $m/k$ is much larger than the typical distance between 
the vertical walls $W(i, k)$
then we will be able to pass through the horizontal stripe at height $n$
only at some exceptional horizontal positions: those places where the
distance between the walls $W(i, k)$ is much larger than typical.
We are only saved if these exceptional positions occur significantly more
frequently than the walls $E(n, k, m)$.

It seems that we need not only bad events, but also some good
ones: the holes.
Our proof systematizes the above ideas by
introducing an abstract notion of walls and holes.
We will have walls and holes of many different types.
To each wall type belongs a fitting hole type.
The walls of a given type occur much less frequently than their fitting
holes.
The whole model will be called a \df{mazery} $\cM$ (a system for creating
mazes). 
The original setting is a very simple mazery:
there is just one wall type, a 1, and one hole type, a 0.
Actually, we have two independent mazeries $\cM_{0}$ and $\cM_{1}$: for
the horizontal and for the vertical line.

\begin{figure}
\if\picsize\SlideSize
    \psset{unit=0.25cm}
\else\if\picsize\PaperSize
    \psset{unit=0.4cm}
\fi\fi

  \includegraphics{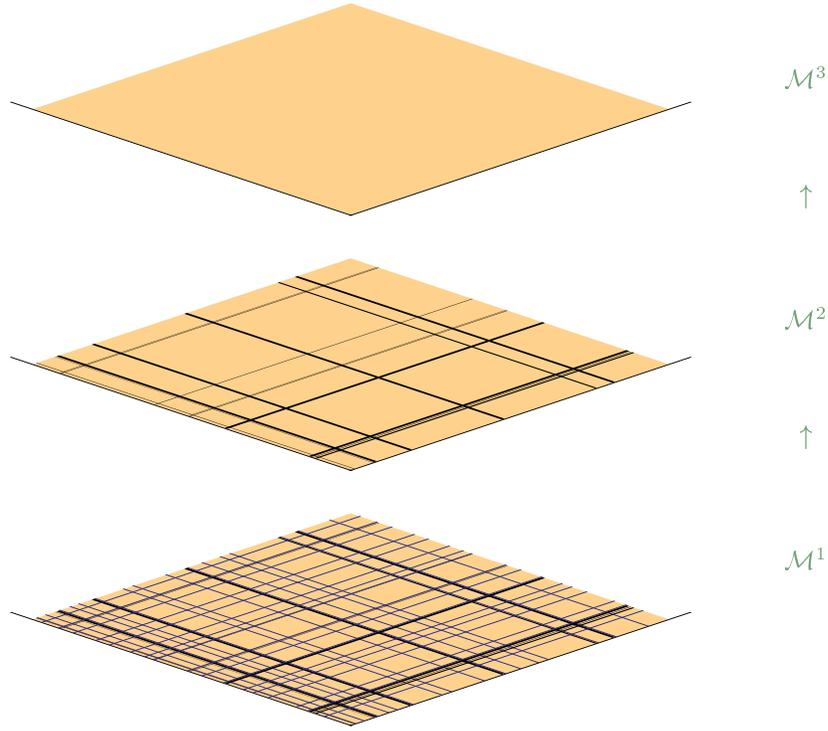}

 \caption{ \label{fig.walls}
Walls in mazeries $\cM_{1},\cM_{2},\cM_{3}$.
The walls seen in $\cM_{2}$ are of the compound type.
(Walls of the emerging type are too rare to show them in a realistic figure.)}
\end{figure}

In any mazery, whenever it happens that
walls are well separated from each other and
holes are not missing, then paths can pass through.
Sometimes, however, unlucky events arise.
It turns out that they can be summarized in two typical examples:
first, when two walls occur too close together (see Figure~\ref{fig.walls});
second, when there is a 
large segment from which a certain hole type is missing.
For any mazery $\cM$, we will define a mazery $\cM^{*}$ whose walls
correspond to these typical unlucky events.
A pair of uncomfortably close walls of $\cM$ gives rise to a wall of
$\cM^{*}$ called a \df{compound wall}.
A large interval without a certain type of hole of $\cM$ gives rise to a
wall of $\cM^{*}$ called an \df{emerging wall}.
Corresponding holes are defined, and it will be shown that the new mazery
also has the property that its holes are much more frequent than the
corresponding walls.
Thus, the ``bad events'' of the outline in
Subsection~\ref{ss.renorm} are the walls of $\cM$, 
the ``very bad events'' are (modulo some details that are not important
now) the compound and emerging walls of $\cM^{*}$.
Let $\cF, \cF'$ be the events $\cF_{k},\cF'_{k}$ formulated in
Subsection~\ref{ss.renorm}.
Thus, $\cF'$ says that in either $\cM_{0}$ or $\cM_{1}$,
a wall appears on the interval $[0, \bub^{*}]$.

The idea of the scale-up construction is that on the level of
$\cM^{*}$ we do not want to see all the details of $\cM$.
We will not see all the walls; however,
some restrictions will be inherited from them: these are distilled in the
concepts of a \df{clean point} and of a \df{slope constraint}.
A point $(x_{0},x_{1})$ is clean for $\cM_{0}\times\cM_{1}$ if $x_{\d}$
is  clean for $\cM_{\d}$ for $\d=0,1$.
Let 
 \[
   \cQ
 \]
be the event that point $(0,0)$ is not clean in $\cM_{0}$ or $\cM_{1}$.

We would like to say that in a mazery, if points $(x_{0},x_{1})$,
$(y_{0},y_{1})$ are such that for $\d=0,1$ we have
$x_{\d} < y_{\d}$ and there are no walls between $x_{\d}$ and $y_{\d}$, 
then $(x_{0},x_{1})\leadsto(y_{0},y_{1})$.
However, this will only hold with some restrictions.
What we will have is the following, with an appropriate parameter
 \[
   0 \le \slb < 1/2.
 \]
 \begin{condition}\label{cnd.reachable-hop}
Suppose that points $(x_{0},x_{1})$, $(y_{0},y_{1})$ are 
such that for $\d=0,1$ we have
$x_{\d} < y_{\d}$ and there are no walls between $x_{\d}$ and $y_{\d}$.
If these points are also clean and satisfy the slope-constraint
 \begin{equation*}
  \slb \le \frac{y_{1}-x_{1}}{y_{0}-x_{0}} \le 1/\slb
 \end{equation*}
then $(x_{0},x_{1})\leadsto(y_{0},y_{1})$.
 \end{condition}

We will also have the following condition:

  \begin{condition}\label{cnd.dense}
Every interval of size $3\bub$ that does not intersect walls
contains a clean point in its middle third.
  \end{condition}

 \begin{lemma}\label{l.imply-reachable}
We have
 \begin{equation}\label{e.imply-reachable}
  \cF \sbs \cF' \cup \cQ.
 \end{equation}
 \end{lemma}
 \begin{proof}
Suppose that $\cQ$ does not hold, then $(0,0)$ is clean.
Suppose also that $\cF'$ does not hold: then by Condition~\ref{cnd.dense}, 
for $\d=0,1$, there is a point $x = (x_{0},x_{1})$ with
$x_{\d} \in [\bub, 2\bub]$ clean in $\cM_{\d}$.
This $x$ also satisfies the slope condition
$1/2 \le x_{1}/x_{0} \le 2$ and is hence, by
Condition~\ref{cnd.reachable-hop}, reachable from $(0,0)$.
 \end{proof}
 
We will define a sequence of mazeries $\cM^{1},\cM^{2},\dotsc$ with
$\cM^{k+1}=(\cM^{k})^{*}$, with $\bub_{k} \to \infty$.
All these mazeries are on a common probability space, since $\cM^{k+1}$ is
a function of $\cM^{k}$.
Actually, we will have two independent mazeries
$\cM^{k}_{\d}$ for $\d=0,1$, constructed from $\cM_{0},\cM_{1}$.
All ingredients of the mazeries will be indexed correspondingly: 
for example, the event $\cQ_{k}$ that $(0,0)$ is not upper right clean in
$\cM_{k}$ plays the role of $\cQ$ for the mazery$\cM^{k}$.
We will have the following property:

 \begin{condition}\label{cnd.Qbd}
 \begin{equation}\label{e.Q-bd}
  \cQ_{k} \sbs \bigcup_{i<k} \cF'_{i}.
 \end{equation}
 \end{condition}
This, along with~\eqref{e.imply-reachable} implies
$\cF_{k} \sbs \bigcup_{i \le k} \cF'_{i}$, which is
inequality~\eqref{e.bd-non-reach-by-wall-all}.
Hence the theorem is implied by the following lemma, which will be
proved after all the details are given:

 \begin{figure}

  \includegraphics{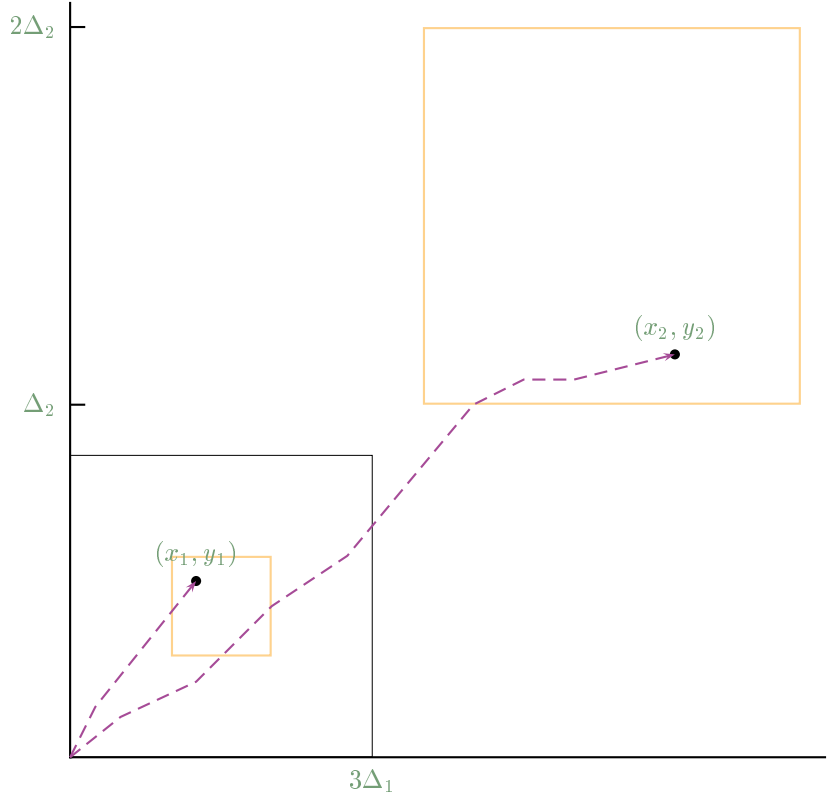}
 
\caption{ \label{fig.main-lemma}
Proof of Theorem~\protect\ref{t.main} from Lemma~\protect\ref{l.main}}
 \end{figure}

 \begin{lemma}[Main]\label{l.main}
If $\p$ is sufficiently small then the sequence $\cM^{k}$ can be
constructed, in such a way that it satisfies all the above conditions and
also 
 \begin{equation}\label{e.main}
   \sum_{k} \Prob(\cF'_{k}) < 1.
 \end{equation}
 \end{lemma}

Figure~\ref{fig.main-lemma} illustrates the proof of the theorem from the
main lemma.

 \begin{remark}
The first clairvoyant synchronization problem, also posed by Winkler,
is more difficult than this one: 
it will be discussed in Section~\ref{s.related}.
 \end{remark}

\subsection{How much lookahead is needed?}\label{ss.lookahead}

Theorem~\ref{t.main} says that when $\p$ is small then with positive
probability, the sequences
$X,Y$ are compatible: they can be synchronized by a clairvoyant demon.
It is a natural question to ask how far a lookahead is needed?
Suppose that sequences $X,Y$ are given up to members $X(n),Y(n)$, and
$X=0$.
Up to which members $X(f(n)), Y(f(n))$ need the sequences to be known in
order to know whether to cast out member $X(n)$?
The following modification of the previous argument gives an upper bound
$f(n)\approx n^{1.5}$.

For point $u \in [0, \bub_{k}]^{2}$, let 
 \[
  \cF_{k}(u)
 \]
be the event that point $u$ is blocked in $[0, \bub_{k}]^{2}$.
Now we assume 
 \[
  \bub_{k+1} > 4\bub_{k},
 \]
and let 
 \[
   \cF''_{k}
 \]
hold iff some wall of $\cM^{k}$ or $\cM^{k+1}$ appears
in $[0, 4 \bub_{k+1}]^{2}$.
Instead of Lemma~\ref{l.imply-reachable}, we have the following.

 \begin{lemma}\label{l.imply-reachable-new}
If $\cF''_{k}$ does not hold then for each
$u \in [0, \bub_{k}]^{2}$ clean in $\cM^{k}$ there is a
$v \in [2\bub_{k}, 3\bub_{k}]^{2}$ clean in $\cM^{k+1}$ and
reachable from $u$.
 \end{lemma}
 \begin{proof}
By Condition~\ref{cnd.dense}, 
for $\d=0,1$, there is a point $v = (v_{0},v_{1})$ with
$v_{\d} \in [2\bub_{k+1}, 3\bub_{k+1}]$ clean in $\cM_{\d}^{k+1}$.
This $v$ also satisfies the slope condition
$1/2 \le (v_{1}-u_{1})/(v_{0}-u_{0}) \le 2$ and is hence, by
Condition~\ref{cnd.reachable-hop}, reachable from $u$.
 \end{proof}

The statement $\sum_{k} \Prob(\cF''_{k}) < 1$ will clearly be proved in just
the same way as the statement $\sum_{k}\Prob(\cF'_{k}) < 1$.
So, with positive probability, none of the events $\cF''_{k}$ occurs.
Then, the above lemma allows us to find a sequence of points
$(0,0) = u_{1}, u_{2}, \dotsc$ such that 
$u_{k} \sbs [2\bub_{k}, 3\bub_{k}]^{2}$, and $u_{k}\leadsto u_{k+1}$.
Finding a path from $u_{k} = (x_{k}, y_{k})$ reaching $u_{k+1}$
means deciding which elements $X(i),Y(j)$ ($x_{k} \le i < x_{k+1}$, 
$y_{k} \le j < y_{k+1}$)
of the sequences $X,Y$ should be cast out.
For this, the sequences $X,Y$ need to be known up to $4\bub_{k+1}$.
Thus, in order to decide about an index $i \ge 2\bub_{k}$, we may need to
know members with indexes $< 4\bub_{k+1}$.
Later in the paper, we will set $\bub_{k+1} = \bub_{k}^{1.5}$.
This shows that lookahead up to index $2 i^{1.5}$ is sufficient.
Different choices of some parameters would allow us to decrease this to
$i^{1+\eps}$ for any $\eps > 0$.

\subsection{The rest of the paper}

In Section~\ref{s.tech-diff}, we discuss the reasons for the less
obvious features of the construction that follows.
The rest of the paper does not depend on this section, it can be skipped,
but it may be helpful to refer back to it for better understanding.

In Section~\ref{s.mazery}, mazeries will be defined.
It is tough to read this section before seeing the role of each notion and
assumption.
The reader may want to refer forward to Section~\ref{s.plan} for such
insights.
 
Section~\ref{s.plan} 
handles all the combinatorial details separable from the calculations.

Section~\ref{s.bounds} defines emerging and compound walls and holes
precisely and estimates their probabilities as much as possible without
bringing in dependence on $k$.

Section~\ref{s.params} defines the scale-up functions,
substitutes them into the earlier estimates and 
finishes the proof of Lemma~\ref{l.main}.

Section~\ref{s.grate-proof} gives a proof left from Section~\ref{s.plan}.

Section~\ref{s.related} relates the present problem to an earlier defined 
clairvoyant synchronization problem.

Section~\ref{s.concl} discusses possible sharpenings of the result.

\section{Some technical difficulties and their solution}\label{s.tech-diff}

In this section, we discuss the reasons for the less
obvious features of the construction that follows.
The rest of the paper does not depend on this section, it can be skipped,
but it may be helpful to refer back to, for better understanding.

\subsection{Clean points}\label{ss.clean}

Our main tool for estimating probabilities will be to attribute various
events to disjoint open intervals.
Cleanness is a property of a point: in order to fit it into this scheme,
it will be broken up into two properties: \df{left-cleanness} and
\df{right-cleanness}.
A point is clean if it is both left-clean and right-clean.
An interval is \df{inner-clean} if its left endpoint is right-clean and its
right endpoint is left-clean; it is outer-clean if its left end is
left-clean and its right end is right-clean.
An inner-clean interval containing no walls will be called a \df{hop}.
(Actually, we will use a slight variant of this notion, called a ``jump'',
see below.)

\subsection{Overlapping walls}

In Section~\ref{s.outline}, we said that when two (vertical) walls
$W_{i}=(x_{i},x_{i}+w_{i})$ ($i=1,2$) 
of mazery $\cM^{k}_{0}$ are too close, they will give rise to a
wall of $\cM^{k+1}_{0}$ of a compound type: 
this postpones the difficulty of dealing with this
situation to a higher level.
It is expected that we can get through these two walls wherever two
fitting (horizontal) holes,
$H_{i}=(y_{i}, y_{i}+h_{i})$ in $\cM^{k}_{1}$, occur at a comparable
distance to each other.
There is much vagueness in these ideas.
What does ``comparable distance'' mean?
What does ``to be expected'' mean?

Let us deal first with the issue: why is it to be expected?
It is true that by the assumptions, $H_{1}$ allows us to get from
$(x_{1},y_{1})$, to $(x_{1}+w_{1}, y_{1}+h_{1})$,
but how do we get from there to $(x_{2}, y_{2})$?
Why can we even assume that walls $W_{1},W_{2}$ are disjoint?
Indeed, suppose that there were three close walls $V_{1},V_{2},V_{3}$ on
level $\cM^{k-1}$, then they would give rise to the compound walls
$V_{1}+V_{2}$ and $V_{2}+V_{3}$, which are not disjoint.

We will deal with the issue by a method used in
the paper repeatedly: we ``define it away''.
It does not exist on level 1.
We require it to be solvable on level $k$ and then use induction to
show that it remains solvable on level $k+1$.
For compound holes, we simply include a requirement into the definition of
the compound hole $H_{1}+H_{2}$, that the interval between $H_{1}$ and
$H_{2}$ is a hop.
This will make it harder to lowerbound the probability of compound holes;
see later.

We require that every interval covered by walls can also be covered by
an interval spanned by a sequence of disjoint walls separated by hops.
In order to prove the same property for $\cM^{k+1}$, it is necessary to
introduce \df{triple compound walls} (essentially, the compounding
operation will be performed twice).
Thus, a sequence of close walls $V_{1},\ldots,V_{n}$ on level $k$ can be
subdivided into a number of disjoint neighbor compounds.
If $n$ is even then all these compounds are of the form $V_{i}+V_{i+1}$.
If $n$ is odd then one of them has the form $V_{i}+V_{i+1}+V_{i+2}$.

This takes care of compound walls, but there are also walls of the emerging
type.
How can these be assumed to be disjoint from everything else and separated
by hops from them?
Again, we define the problem away: we will essentially introduce emerging
walls one-by-one, allowing one only if it is disjoint from the others
and is outer-clean.

The last step breaks an important property of our model.
It is a global operation: now
whether an interval $I$ is the body of a wall of $\cM^{k}_{0}$
is no more an
event depending only on the random variables $X_{0}(j)$ in this interval.
The events that there is a wall of a certain type on interval
$I_{1}$ and of some other type on interval $I_{2}$, are not independent
anymore: we seem to lose our only tool of probability estimation.
Fortunately, this problem can also be defined away, since
for walls, we need only probability upper bounds.
We introduce another concept, the concept of a \df{barrier}.
Every wall is a barrier, but not vice versa.
The occurrence of a barrier of any  given 
type on an open interval $I$ will only depend on
elementary events in $I$, and will be independent of the occurrence of a
barrier of any other type on any interval $J$ disjoint from $I$.
Our probability bounds will be for barriers.

\subsection{Compound wall types}

When we say that walls are much less probable than the fitting holes then
what we mean is that there is an overall constant 
 \begin{equation}\label{e.hxp}
  \hxp = 0.03
 \end{equation}
such that when
$p$ is the probability of a certain wall type and $q$ the probability of a
fitting hole type then $q \ge p^{\hxp}$.

\begin{figure}


 \includegraphics{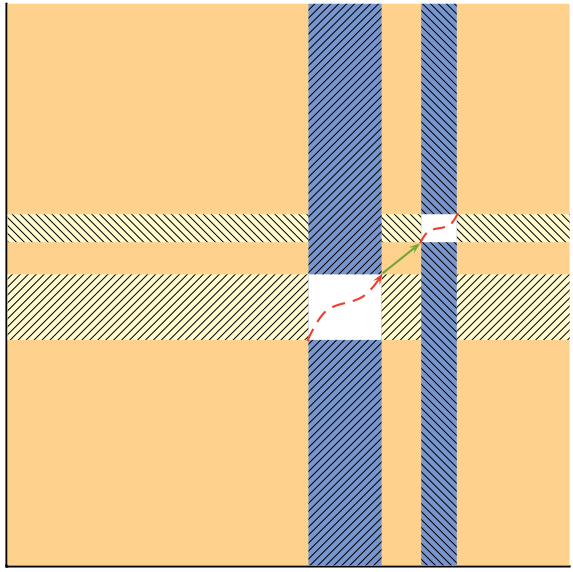}

 \caption{ \label{fig.compound} 
Compound wall and a fitting compound hole}
 \end{figure}

Several times, it will be important to bound the probability that any barrier
at all occurs at a given point.
Since our probability upper bounds apply generally to barriers
of some given type, it is important to limit the number of types.
We said that given two barriers $W_{i}$ $(i=1,2)$ of types $\ag_{i}$
at a certain small distance $d$, a compound barrier should arise.
``Small'' will mean here smaller than a certain parameter $\f$.
We cannot afford a new compound barrier type for each triple
$(\ag_{1},\ag_{2},d)$ for all $d < \f$: this would lead to too many types.
On the other hand, if we only introduce a single type for all pairs
$\ag_{1},\ag_{2}$ then a simple probability upper bound of this type would
only be $\f \p(\ag_{1})\p(\ag_{2})$, while the probability lower bound on
the fitting compound hole type would only be
$\p(\ag_{1})^{\hxp}\p(\ag_{2})^{\hxp}$, ignoring the factor $\f$.
We need a compromise between these two extremes.
The parameter 
 \begin{equation}\label{e.lambda-def}
  \lg = 2^{1/4}
 \end{equation}
 will be introduced, and we will have compound barrier types 
 \begin{equation}\label{e.compound-hole-i}
 \bg = \ang{\ag_{1},\ag_{2},i}
 \end{equation}
for different $i$.
A pair of barriers $W_{i}$ of type $\ag_{i}$
gives rise to a compound barrier of such a type if their distance falls
between $\lg^{i}$ and $\lg^{i+1}$.
A fitting compound hole will have two fitting holes,
whose distance is between the values $\lg^{i-1}$ and $\lg^{i}$.
(See Figure~\ref{fig.compound}.)

With the above choice of barrier types, even barriers belonging to
a given type can have different sizes.
This makes it harder to prove even a probability upper bound
$\p(\ag_{1}) \p(\ag_{2})$ on the occurrence of barriers $W_{i}$ of type
$\ag_{i}$ adjacent to each other.
We again define the problem away:
the probability bound $\p(\ag)$ on the occurrence of a barrier of type $\ag$
is simply defined as the sum of bounds $\p(\ag,w)$ of the occurrence of
type $\ag$ and size $w$.

\subsection{Ranks}

Having all these different compound barrier types implies that 
even if the component types have probability $\p(\ag_{i}) = s$, the probability
bound on compound barriers now ranges from $s^{2}$ to $\f s^{2}$.
Assuming that $\f=s^{-\fxp}$ for some small $\fxp$, the probability bound
ranges from $s^{2}$ to $s^{2 - \fxp}$.
Repeating this operation $k$ times, the bound ranges from $s^{2^{k}}$ to 
$s^{(2 - \fxp)^{k}}$.
This makes the smallest probability of a barrier type much smaller than the
probability of all barrier types.

Recall that a (vertical) barrier of the emerging type
occurs in $\cM^{k+1}_{0}$ on
an interval $I$ of size $\g$ (for a certain parameter $\g$ smaller than the
$\f$ introduced above), if there is a 
barrier type $\ag_{0}$ of $\cM^{k}$ such that no hole of fitting type
$\ag'_{0}$ appears in $I$.
A (horizontal) hole of a fitting type of $\cM^{k}_{1}$ 
occurs on an interval $J$ of size $\approx \g$ if
no (horizontal) barrier of any type appears on $J$.
It is not sufficient to exclude barriers of type $\ag_{0}$, since vertical
holes of other types may also be missing from $I$.
A simple probability upper bound for any barrier appearing on $J$ is
$\g \ol \p = \g \sum_{\ag} \p(\ag)$.
This can only be small if $\g < 1 / \ol \p$.
The probability that $I$ is an emerging barrier is exponentially small in $\g$
but this becomes significant only if $\g$ is larger than the inverse of the
probability lower bound $(\p(\ag_{0}))^{\hxp}$ on the appearance of a hole of
type $\ag'_{0}$. 
Thus, we need 
 \begin{equation}\label{e.impossible}
 \ol \p < (\p(\ag_{0}))^{\hxp}, 
 \end{equation}
which is impossible
if the difference between the smallest and the largest $\p(\ag)$ grows as
shown above.

The solution is to realize that if we have barrier types with many different
probability bounds, we do not have to deal with them all at the same time.
We assign a number called \df{rank} to each barrier type, based essentially
on its probability bound.
The probability bound of barriers of rank $\r$ will be 
$\p(\r) \approx \lg^{-\r}$, and this will also bound approximately the
probability that any barrier of rank $\r$ starts at a given point.
Mazery $\cM^{k}$ will have types present in it that are between ranks
$\R_{k}$ and $\R_{k}^{c}$ for a certain constant $c$, where
$\R_{k}=\lg^{\txp^{k}}$, $\txp=2 - \fxp$.
When going from level $k$ to level $k+1$, we eliminate only the ranks
between $\R_{k}$ and $\R_{k+1}$ (and we add some new, higher ranks).
The barrier type $\ag_{0}$ above in the introduction of an emerging
barrier, will be restricted to have rank less than $\R_{k+1}$ (such barrier
types will be called \df{light}).
With the rank limited from above, the fitting hole
probability lower bound is limited from below, and we escape the
impossible requirement~\eqref{e.impossible}.

\subsection{Compound hole probability}

It is relatively easy to obtain an upper bound
$(\lg^{i+1} - \lg^{i})\p(\ag_{1})\p(\ag_{2})$ 
on the probability of a compound barrier
of type $\bg$ in~\eqref{e.compound-hole-i}.
Indeed, we essentially sum up the probability bound
$\p(\ag_{1})\p(\ag_{2})$ for the different possible values of the distance
$d \in [\lg^{i}, \lg^{i+1}-1]$ between the two barriers.
It is less clear how to find the corresponding lower bound
 \begin{equation}\label{e.compound-hole-inform}
  \ul p = (\lg^{i+1}-\lg^{i})^{\hxp}(\p(\ag_{1}))^{\hxp}(\p(\ag_{2}))^{\hxp} 
 \end{equation}
for the fitting hole.
In the upper bound, we can ignore the fact that the various possibilities
for distance $d$ are not disjoint events; not so for the lower bound.
We will use a combination of two methods.
If $i$ is small then we will use the main method of the paper:
define the problem away, by requiring something in the definition of
a mazery that will imply the result (see the hole lower
bound~\eqref{e.hole-lb}).
If $i$ is large, then we break up the interval considered into
disjoint subintervals of size $3\bub$, and use the independence of the relevant
events on these along with some computation, in order to infer the same
property~\eqref{e.hole-lb} for $\cM^{k+1}$.

The above tactic gives the desired lower bound on the probability that two
holes of the given types appear at a distance 
$d \in [\lg^{i-1}, \lg^{i} - 1]$, but
it still does not guarantee that the two holes are separated by an
interval free of barriers (inner-cleanness of this interval is easier to
achieve). 
The probability that an interval of size $\f$ 
is barrier-free, has a large lower bound: the
problem is only to be able to multiply this lower bound with the lower
bound $\ul p$ in~\eqref{e.compound-hole-inform}.
We will make use of the FKG inequality for
this, noticing that the events whose intersection we want to form are
decreasing functions of the original sequence $X_{0}$.
(It is possible to
avoid the reference to monotonicity, at the price of a little more sweat.)

\section{Walls and holes}\label{s.mazery}

 \subsection{Notation}
  
We will use
 \[
   a \et b = \min(a,b),\quad a \vee b = \max(a,b).
 \]
The \df{size} of an interval $I = (a,b)$ is denoted by $|I| = b-a$.
For two sets $A, B$ in the plane or on the line,
 \[
   A + B = \setof{a + b : a \in A,\; b \in B}.
 \]
For two
different points $u_{i} = (x_{i},y_{i})$ $(i=0,1)$ in the plane, when 
$x_{0} \le x_{1}$, $y_{0} \le y_{1}$:
 \begin{align*}
      \slope(u_{0}, u_{1}) &= \frac{y_{1}-y_{0}}{x_{1}-x_{0}}, 
\\ \minslope(u_{0}, u_{1}) &= 
   \min\bigparen{\slope(u_{0},u_{1}), 1/\slope(u_{0},u_{1})}.
  \end{align*}

\subsection{Mazeries}

\subsubsection{The whole structure}

A \df{mazery} 
 \begin{equation}\label{e.mazery}
 \MM = (\Btypes, \Htypes, \Rank(\cdot), \Bvalues, \Hvalues, \cM, \bub,
\slb, \p(\cdot)	) 
 \end{equation}
has the following parts.
 \begin{enumerate}[--]

  \item Two disjoint finite sets $\Btypes$ and $\Htypes$ called the sets of
\df{barrier types} and \df{hole types}.

  \item A function $\Rank: (\Btypes \cup \Htypes) \to \ZZ_{+}$.
Each type $\ag$ has a \df{rank} $\Rank(\ag)$ that is a nonnegative
integer.

  \item 
Let $\cI$ be the set of nonempty
open intervals with integer endpoints, then we have two sets:
 \[
      \Bvalues \sbs \cI \times \Btypes,
\quad \Hvalues \sbs \cI \times \Htypes.
 \]
Elements of $\Bvalues$ and $\Hvalues$ are called \df{barrier values} and
\df{hole values} respectively.
A barrier or hole value $E = (B, \ag)$
has \df{body} $B$ which is an open interval, and \df{type} $\ag$.
We write $\Body(E) = B$, $|E| = |B|$.
We will sometimes denote the body also by $E$.
It is not considered empty even if it has size 1.
We write $\Type(E) = \ag$.

  \item The random process $\cM$ will be detailed below.

  \item The parameters $\bub >0$, $\slb \ge 0$ and the 
probability bounds $\p(\r)$ will also be detailed below.

 \end{enumerate}

The following conditions hold for the parts discussed above.

 \begin{condition}\label{cnd.types}\
 \begin{enumerate}[\upshape 1.]

  \item\label{i.types.bub}
Let us denote by $\sizelb(\ag),\sizeub(\ag)$ the infimum and
supremum of the sizes of barriers or holes of type $\ag$.
We require
 \[
   \sizeub(\ag) \le \bub.
 \]

  \item\label{i.types.fitting}
To each barrier type $\ag$ there corresponds a \df{fitting} hole type $\ag'$
with
 \[
   \sizeub(\ag') \le \sizelb(\ag), \quad \Rank(\ag') = \Rank(\ag). 
 \]
(We require a hole to have a smaller or equal size than the barrier it
fits; this will bound the amount by which the minimal slope of a path 
must increase while passing a wall.  
Equality is achieved in the example
where a wall of $k$ 1's can be passed only at a hole of $k$ 0's: thus, this
example is the worst that can happen.)

 \end{enumerate}
 \end{condition}

\subsubsection{The random processes}

Now we discuss the various parts of the random process
 \begin{align*}
   \cM  &= (Z, \cB, \cW, \cH, \cC, \cS).
 \end{align*}
Here, $Z = (Z(0), Z(1), \dotsc)$, is 
a sequence of independent, equally distributed 0-1-valued random
variables. 
Since $Z(i)$ is thought of as belonging to the open interval $(i,i+1)$, we
will write for any interval $(a, b)$:
 \[
   Z((a, b)) = (Z(a), \ldots, Z(b-1)).
 \]
We will have the following random sets:
 \begin{equation*}
          \cW \sbs \cB   \sbs \Bvalues,
\quad              \cH   \sbs \Hvalues,
\quad     \cS \sbs \cC   \sbs \ZZ_{+} \times (\ZZ \xcpt \{0\}),
 \end{equation*}
all of which are functions of $Z$.
The elements of $\cB$ are \df{barriers}.
The elements of $\cW$ are called \df{walls}: each wall is also a barrier.
The elements of $\cH$ are called \df{holes}.

 \begin{remarks}\label{r.wall-barrier}\
  \begin{enumerate}[\upshape 1.]
   \item 
In what follows we will refer to $\cM$ by itself 
also as a \df{mazery}, and will mention $\MM$ only rarely.
This should not cause confusion; though $\cM$ is part of $\MM$,
it relies implicitly on all the other ingredients of $\MM$.
 
    \item
The distribution of $\cB$ is simpler than that of $\cW$, 
but sample sets of $\cW$ will have a more useful structure.
The parts of the paper dealing with combinatorial questions (reachability)
will work mainly with walls and hops (see below),
the parts containing probability calculations will
work with barriers and jumps (see below).

  \end{enumerate}
 \end{remarks}

\subsubsection{Cleanness}

For $x \ge 0, r > 0$,
if $(x,-r) \in \cC$ then we say that point $x$ is \df{left $r$-clean}.
If $(x,r) \in \cC$ then we say that $x$ is \df{right $r$-clean}.
A point $x$ is called left-clean (right-clean) 
if it is left $r$-clean (right $r$-clean) for all $r$.
It is \df{clean} if it is both left- and right-clean.
If the left end of an interval $I$ is right $|I|$-clean and
the right end is left $|I|$-clean then we say $I$ is \df{inner-clean}.
If its left end is left-clean and its right end is right-clean then we say
that it is \df{outer-clean}.
Let $(a,b), I$ be two intervals with $(a-\bub, b+\bub) \sbs I$.
We say that $(a,b)$ is \df{cleanly} contained in $I$ if it is outer-clean.

For $x \ge 0, r > 0$, if $(x,-r) \in \cS$ then we say that point $x$ is 
\df{strongly} left $r$-clean.
Similarly, every concept concerning cleanness has a counterpart when we 
replace $\cC$ with $\cS$ and clean with strongly clean.

A closed interval is called a \df{hop} if it is inner clean and contains no
wall.
It is a \df{jump} if it is strongly inner clean and contains no barrier.
By definition, all jumps are hops.
A hop or jump may consist of a single point.

Let us call an interval \df{external} if it does not intersect any wall.
(Remember that the bodies of walls and holes are open intervals.)
Two disjoint walls or holes are called \df{neighbors} if the interval
between them is a hop.
A sequence $W_{i} \in \cW$ of walls
$i=1,2,\dotsc$ is called a \df{sequence of neighbor walls}
if for all $i$, $W_{i}$ is a left neighbor of $W_{i+1}$.

\subsubsection{Conditions on the process $\cM$}

There are many conditions on the distribution of process $\cM$, 
but most of them are fairly natural.
We would like to call
attention to a crucial condition that is not clearly motivated:
Condition~\ref{cnd.distr}.\ref{i.distr.hole-lb}.
It implies, as a special case with $c=b+1$,
that through every wall, at every position,
there is a fitting hole with sufficiently large probability.
The general case has been formulated carefully and it
is crucial for the inductive proof that the hole lower bound will also hold
on compound walls after renormalization (going from $\cM^{k}$ to
$\cM^{k+1}$). 

The function 
 \begin{equation}\label{e.p(,)}
  \p(\r,w)
 \end{equation}
is defined as the supremum of probabilities (over all points)
that any barrier with rank $\r$ and size $w$ starts at a given point.
The function $\p(\r)$ will be an upper bound on $\sum_{w} \p(\r, w)$.
Let
 \[
   \p(\ag, w)=\p(\r, w),\quad \p(\ag) = \p(\r)   
 \]
for any type $\ag$ of rank $\r$.

Let us introduce the constants
 \begin{equation}\label{e.Rlb.aux1}
  \aux_{1} = 6, \quad \aux_{0} \approx 7.41,
 \end{equation}
where the requirement $\aux_{1} > 5$ comes from 
inequality~\eqref{e.Rlb.compound-contrib} below, while
the value for $\aux_{0}$ will be motivated in~\eqref{e.Rlb.aux3}.
For each rank $\r$, let us define the function
 \begin{equation}\label{e.h-def}
            \h(\r)  = \aux_{0}\r^{\aux_{1}\hxp}(\p(\r))^{\hxp},
\quad       \h(\ag) = \h(\Rank(\ag)).
 \end{equation}
The exponent $\hxp$ has been introduced in~\eqref{e.hxp}.
Its choice will be motivated in Section~\ref{s.params}.
The factor $\aux_{0}\r^{\aux_{1}\hxp}$ will absorb some 
nuisance terms as they arise in the estimates.
The function $\h(\r)$ will be used as a lower bound for the probability of
holes of rank $\r$.

 \begin{condition}\label{cnd.distr}\
 \begin{enumerate}[\upshape 1.]

  \item\label{i.distr.indep}(Dependencies and monotonicity)
  \begin{enumerate}[\upshape a.]

    \item\label{i.distr.indep.barr}
For a barrier value $E$, the event $\setof{E \in \cB}$
is an increasing function of $Z(\Body(E))$.
Thus, any set of events of the form $E_{i} \in \cB$
where $E_{i}$ are disjoint, is independent.

    \item\label{i.distr.indep.hole}
For a hole value $E$, the event $\setof{E \in \cH}$
is a decreasing function of $Z(\Body(E))$.
   
    \item\label{i.distr.indep.clean}
For every point $x$ and integer $r$, the events
$\setof{(x,-r) \in \cS}$, $\setof{(x, r) \in \cS}$  
are decreasing functions of $Z((x-r,x))$ and $Z((x,x+r))$ respectively.

When $Z$ is fixed, strong and not strong left (right) $r$-cleanness are
decreasing as functions of $r$.
These functions reach their minimum at $r=\bub$: thus, if 
$x$ is (strongly) left (right) $\bub$-clean then it is (strongly)
left (right)-clean.
  \end{enumerate}

  \item\label{i.distr.combinat}(Combinatorial requirements)
   \begin{enumerate}[\upshape a.]

    \item\label{i.distr.inner-clean}
A maximal external interval is inner-clean (and hence is a hop).

    \item\label{i.distr.cover}
If an interval $I$ is surrounded by maximal external intervals of size 
$\ge \bub$ then it is spanned by a sequence of neighbor walls.
(Thus, it is covered by the union of the neighbor walls of this sequence
and the hops between them.)

  \item\label{i.distr.clean}
If an interval of size $\ge 3\bub$ contains no walls then
its middle third contains a clean point.

   \end{enumerate}

  \item\label{i.distr.bounds}(Probability bounds)

  \begin{enumerate}[\upshape a.]

   \item  For all $\r$ we have 
 \begin{equation}\label{e.p()}
  \p(\r) \ge \sum_{l} \p(\r,l).
 \end{equation}

   \item
The following requirement imposes an implicit upper bound on $\p(\ag)$:
  \begin{equation}\label{e.0.6}
   \bub^{\hxp} \h(\ag) < 0.6.
  \end{equation}

   \item\label{i.distr.bounds.ncln-ub}
For $\ncln = \sup_{x}\Prob\setof{x \txt{ is not strongly clean}}$, we have
  \begin{equation}\label{e.ncln-ub}
  \ncln < 0.25.
  \end{equation}

   \item\label{i.distr.hole-lb}
Let $\ag$ be a barrier type, let $a \le b < c$ and $b-a, c-b  \le 6 \bub$,
and let $E(a, b, c, \ag)$ be the event
that there is a $d \in [b, c-1]$ such that $[a,d]$ is a jump
and a hole of type $\ag'$ starts at $d$.
Then
 \begin{equation}\label{e.hole-lb}
   \Prob(E(a,b,c,\ag)) \ge (c-b)^{\hxp} \h(\ag).
 \end{equation}
  
 \end{enumerate}

 \end{enumerate}
 \end{condition}

\begin{remark}\label{r.distr}
These conditions imply that the property that an interval $I$ is strongly
cleanly contained in some interval $J$ depends only on $Z(J)$.
\end{remark}

\subsubsection{Reachability}
  
Take two random processes
 \[
  \cM_{\d} \quad (\d=0,1),
 \]
independent, and distributed like $\cM$.
To the pair of mazeries $\cM_{0}, \cM_{1}$ belongs a random graph 
 \begin{equation*}
       \cV = \ZZ_{+}^{2},
\quad  \cG = (\cV, \cE)
 \end{equation*}
where $\cE$ is determined by the  above random processes as in 
Subsection~\ref{ss.perc}.
From now on, \df{reachability} is always understood in the graph $\cG$.
Just as the random sets $\cW, \cH$ of walls and holes were introduced as
parts of what makes up a mazery $\cM$, for example the random set
$\cW_{0}$ is the corresponding part of mazery $\cM_{0}$.
If $u = (x_{0}, x_{1})$ and $v = (y_{0}, y_{1})$ 
are points of $\cV$ such that
for $\d = 0,1$, $x_{\d} < y_{d}$, and $(x_{\d}, y_{\d})$ is a hop
then we will say that $(u, v)$ is a \df{hop}.

Let the interval $I = (a_{1},a_{2})$ be the body of
a (vertical) barrier $B$.
For an interval $J = (b_{1}, b_{2})$ with 
$|J| \le |I|$  we say that $J$
is a (horizontal) \df{hole passing through} $B$, or \df{fitting} $B$, if
$(a_{1}, b_{1})\leadsto(a_{2}, b_{2})$.
This hole is called \df{left clean}, right clean or outer clean
if $J$ has these properties in $\cM_{1}$.
Vertical holes are defined similarly.

The graph $\cG$ is required to satisfy the following conditions.

 \begin{condition}[Reachability]\label{cnd.reachable}\
 \begin{enumerate}[\upshape 1.]

  \item\label{i.reachable.wall}
Let $u = (x_{0}, x_{1})$.
Suppose that a wall $W$ of $\cW_{0}$ starts at
$x_{0}$, and a fitting hole $H$ of $\cH_{1}$ starts at $x_{1}$.
Then $u\; \leadsto\; u + (|W|, |H|)$.
The same holds when we interchange the indexes 0 and 1.

  \item\label{i.reachable.hop}
If $u,v$ are points of $\cV$ such that $(u, v)$ is a hop
and $\minslope(u, v) \ge \slb$, then $u\leadsto v$.
We require
 \begin{equation}\label{e.slb-ub}
  0 \le \slb < 0.5.
 \end{equation}
 \end{enumerate}
 \end{condition}

 \begin{example}\label{x.base}
The compatible sequences
problem can be seen as a special case of such a mazery.
We define this mazery as $\cM^{1}$, the first one of a series of mazeries
$\cM^{k}$ to be defined later.
There is only one barrier type, of rank $\R_{1}$ where the number $\R_{1}$
will be chosen conveniently (sufficiently large) later, and one hole type.
We have $\Delta =1$ and $\slb = 0$.
The reader should check that if
the probability of the barrier type $\p$ is chosen sufficiently small
then the hole lower bound~\eqref{e.hole-lb} and the bound~\eqref{e.0.6}
will still be satisfied.
There is a wall of size 1 starting at $j$ if $Z(j) = 1$, and a hole of size
1 if $Z(j) = 0$.
Barriers are walls, and every point is strongly clean.
 \end{example}

 \section{The scaled-up structure}\label{s.plan}

In this section, we will define the structural parts of the
scaling-up operation $\MM \mapsto \MM^{*}$: 
we still postpone the
definition of various parameters and probability bounds for $\MM^{*}$.
We start by specifying when we will consider walls to be sufficiently
sparse and holes sufficiently dense, and what follows from this. 
Let $\slopeincr$ be a constant and $\f,\g$ be some parameters with 
 \begin{align}
\label{e.slopeincr}
   \slopeincr                  &= 48,
\\\label{e.g-intro}
                            \g &> 6 \bub,  
\\\label{e.f-intro} \slopeincr \g / \f &< 0.5 - \slb.
 \end{align}
Let
 \begin{equation}\label{e.new-slb}
   \slb^{*} = \slb + \slopeincr \g / \f.
 \end{equation}
We will say that the process $(\cM_{0},\cM_{1})$ satisfies the 
\df{grate condition} over the rectangle $I_{0} \times I_{1}$
with parameters $\f, \g$ if the following holds.

 \begin{condition}[Grate]\label{cnd.grate}\

 \begin{enumerate}[\upshape 1.]

  \item\label{i.grate.pre-hop}
For $\d=0,1$, for some $n_{\d} \ge 0$, there
is a sequence of neighbor walls $W_{\d,1},\ldots,W_{\d,n_{\d}}$
and hops
around them such that taking the walls along with the hops
between and around them, the union is $I_{\d}$
(remember that a hop is a closed interval).
Also, the hops 
between the walls have size $\ge \f$, and the hops
next to the ends of $I_{\d}$ (if $n_{\d} > 0$) have size $\ge \f/3$.

  \item\label{i.grate.hole}
For every wall $W$ occurring in $I_{\d}$,
for every subinterval $J$ of size $\g$ of the hops
between and around
walls of $\cW_{1-\d}$ there is a hole of $\cH_{1-\d}$ fitting $W$, with
its body cleanly contained in $J$.

 \end{enumerate}
 \end{condition}

 \begin{lemma}[Grate]\label{l.grate}
Recall the constant $\slopeincr$ from~\eqref{e.slopeincr}.
Suppose that $(\cM_{0},\cM_{1})$ satisfies the grate condition
over the rectangle $I_{0} \times I_{1}$ with 
$I_{\d} = [u_{\d}, v_{\d}]$, $u = (u_{0}, u_{1})$, $v = (v_{0}, v_{1})$. 
If
 \[
   \minslope(u, v) \ge \slb^{*}
 \]
then $u\leadsto v$ in $\cG$.
 \end{lemma}

We postpone the proof of this lemma to Section~\ref{s.grate-proof}.

After defining the mazery $\MM^{*}$, eventually we will have to prove the
required properties.
To prove Condition~\ref{cnd.reachable}.\ref{i.reachable.hop}
for $\MM^{*}$, we will invoke
the Grate Lemma~\ref{l.grate}: for that, we will need the Grate
Condition~\ref{cnd.grate}, with appropriate parameters $\f,\g$.
To ensure these conditions in $\cM^{*}$, we will introduce some new barrier
types whenever they would fail; this way, they will be guaranteed to hold
in a wall-free interval of $\cM^{*}$.
The set $\Btypes^{*}$ will thus contain a new, so-called 
\df{emerging} barrier type, and several new, so-called \df{compound}
barrier types, as defined below, and 
$\Htypes^{*}$ will contain the corresponding hole types.
The following algorithm creates the barriers, walls and holes of these new
types out of $\cM$.
Let us denote by $\R$ a lower bound on all possible ranks,
and by $\R^{*} > \R$ the number that will become the lower bound in
$\cM^{*}$.
We will also make use of parameter $\lg$ defined in~\eqref{e.lambda-def}
with the property 
 \begin{equation}\label{e.R-cond}
   \R^{*} \le 2 \R - \log_{\lg} \f.
 \end{equation}
Types of rank lower than $\R^{*}$ are called \df{light}, the other ones
are called \df{heavy}.

 \begin{enumerate}[\upshape 1.]

  \item(Cleanness)\label{i.inform.clean}
For an integer $r>0$, a
point $x$ of $\cM^{*}$ will be called $r$-left-clean if it is
$r$-left-clean in $\cM$ and there is no wall of $\cM$ in
$[x-r, x]$ whose right end is closer to $x$ than $\f/3$.
Right-cleanness is defined similarly.

The point $x$ will be called strongly $r$-left-clean in $\cM^{*}$
if it is strongly $r$-left-clean in $\cM$ and 
there is no barrier of $\cM$ in
$[x-r, x]$ whose right end is closer to $x$ than $\f/3$.
Strong right-cleanness is defined similarly.

  \item(Emerging type)\label{i.inform.emerg}
There is a
new barrier type $\ang{\g}$ in $\Btypes^{*}$, called the \df{emerging type},
with rank $\R^{*}$, so it will be heavy.
Suppose that there are $a \le a' < b' \le b$ with
$a'-a \le \bub$, $b'-a' = \g-4\bub$, $b-b' \le \bub$ and a
light barrier type $\ag$ such that 
no hole of type $\ag'$ is cleanly contained in $[a',b']$. 
Then $((a,b),\ang{\g})$ is a barrier of $\cM^{*}$.

Let us list all barriers of type $\ang{\g}$ in a
(typically infinite) sequence $B_{1},B_{2},\dots$.
We will mark some of these barriers as walls of $\cM^{*}$ as follows.
Let us go through the sequence $B_{1},B_{2},\dots$ in order.
Suppose that we have already decided which of the barriers $B_{i}$, $i<n$
are walls.
Then we mark $B_{n}$ with body $[a,b]$ as a wall if and only if the
following conditions hold:
 \begin{enumerate}[$\bullet$]
  \item $[a,b]$ does not intersect any of
the barriers $B_{i}$, $i<n$ already marked as walls;
  \item $[a,b]$ is a hop of $\cM$, cleanly contained in a hop.
 \end{enumerate}
There is a fitting hole type $\ang{\g}'$.
Any interval of size $\g - 4 \bub$ is a hole of this type in $\cM^{*}$, 
if it is a jump of $\cM$.

  \item(Compound types)\label{i.inform.compound}
We make use of a certain sequence of integers $d_{i}$, which is defined by
the following formula, but only for $d_{i} \le \f$:
 \begin{equation}\label{e.jump-types}
   d_{i} = \begin{cases}
               i &\txt{ if $0 \le i < 17$},
\\ \flo{\lg^{i}} &\txt{ if $i \ge 17$}.
           \end{cases}
 \end{equation}
For any pair $\ag_{1},\ag_{2}$ of barrier types where $\ag_{1}$ is light,
and any $0 \le i$ with $d_{i} \le \f$, we introduce a new type
$\bg = \ang{\ag_{1},\ag_{2},i}$ in $\Btypes^{*}$.
(The type $\ag_{2}$ can be not only a heavy type of $\cM$ but also
the new emerging type $\ang{\g}$ just defined above.)
Such types will be called \df{compound types}.
The rank of this type is defined by 
 \begin{equation}\label{e.compound-rank}
 \Rank(\bg) = \Rank(\ag_{1}) + \Rank(\ag_{2}) - i.
 \end{equation}
It follows from~\eqref{e.R-cond} that these new types are heavy.
A barrier of type $\bg$ occurs in $\cM^{*}$ wherever disjoint barriers
$W_{1}, W_{2}$ of types $\ag_{1}, \ag_{2}$ occur (in this order)
at a distance $d \in [d_{i},d_{i+1}-1]$.
(If $d_{i+1}$ is not defined then $d \in [d_{i},\f-1]$.)
Thus, a shorter distance gives higher rank.
The body of this compound barrier is the union of the bodies of 
$W_{1},W_{2}$ and the interval between them.
We denote the new compound barrier by
 \[
  W_{1} + W_{2}.
 \]
It is also a wall if $W_{1},W_{2}$ are walls of $\cM$ or
$\cM^{*}$ (we have already introduced some walls of $\cM^{*}$, of the
emerging type) and the interval between the $W_{i}$ is a hop of $\cM$.
To make sure this construction always succeed, we require
 \begin{equation}\label{e.f-bub-star}
   2\f + \bub < \bub^{*}.
 \end{equation}
Compound hole types will be defined similarly, using a sequence $e_{i}$,
$i=1,2,\dots$ defined by the following formula, but only if
$d_{i} \le \f$:
 \begin{equation}\label{e.hole-jump-types}
  e_{i} = \begin{cases}
            i       & \txt{ if $0 \le i \le 17$},
\\          d_{i-1} & \txt{ if $i > 17$}.
	  \end{cases}
 \end{equation}
Thus, the distance between the component holes of a compound hole 
of type $\bg' = \ang{\ag'_{1},\ag'_{2}, i}$ varies in the interval
$[e_{i},e_{i+1}-1]$.
This hole type fits the barrier type $\bg = \ang{\ag_{1},\ag_{2},i}$.
A hole of type $\bg'$ occurs in $\cM^{*}$
if holes $H_{1}, H_{2}$ of types $\ag'_{1},\ag'_{2}$ occur,
connected by a jump of $\cM$ of size $d \in [e_{i}, e_{i+1}-1]$.

Now we repeat the whole compounding step, introducing compound types in 
which now $\ag_{2}$ is required to be light.
The type $\ag_{1}$ can be any type introduced until now, also a compound
type introduced in the first compounding step.
So, the walls that will occur as a result of the compounding operation are
of the type $L$-$*$, $*$-$L$, or $L$-$*$-$L$, where $L$ is a light wall of
$\cM$ and $*$ is any wall of $\cM$ or an emerging wall of $\cM^{*}$.

Note that emerging walls were made to be disjoint from each other, but
compound walls were not.

  \item(Finish)\label{i.inform.rm-light}
Remove all light types and all the corresponding barriers, walls and holes.
Moreover, if a heavy wall $W$ of $\cM$ is contained in an outer clean light
wall $W'$ of $\cM$ then with $W'$, remove also $W$. 
(We do not know whether such heavy walls can really occur but their removal
will not hurt.)

The graph $\cG$ does not change in the scale-up: $\cG^{*}=\cG$.
 
 \end{enumerate}

Let us prove some of the required properties of $\MM^{*}$.

 \begin{lemma}\label{l.distr.indep}
The new mazery $\cM^{*}$ satisfies
Condition~\ref{cnd.distr}.\ref{i.distr.indep}.
 \end{lemma}
 \begin{proof}
We will see that all the properties in the condition follow essentially
from the form of our definitions.
Note that when an existential or universal quantifier is applied to
a family of
monotonically increasing (decreasing) events, the result is monotonically
increasing (decreasing) event (as a function of the sequence $Z$).
Indeed, these quantifiers are just the maximum and minimum operations.

Condition~\ref{cnd.distr}.\ref{i.distr.indep.barr} says that
for a barrier value $E$, the event $\setof{E \in \cB}$
is an increasing function of $Z(\Body(E))$.
To check this, consider all possible barriers of $\cM^{*}$.
We have the following kinds:
 \begin{enumerate}[--]

  \item
Heavy barrier values $E$ of $\cM$: all heavy barriers of $\cM$ remained 
barriers of $\cM^{*}$, so monotonicity holds, and
the event still depends only on $Z(\Body(E))$.

  \item
Barrier values $E$ of $\cM^{*}$ of the emerging type: 
such a barrier is defined by the
absence of some holes of $\cM$ on some subintervals of $\Body(E)$.
Since the presence of a hole in $\cM$ is a monotonically decreasing event,
the presence of an emerging barrier is an increasing event
depending only on $Z(\Body(E))$.

  \item
Barrier values $E$ of some compound type.
Such a barrier appears in $\cM^{*}$ if certain barriers appear in $\cM$ in
$\Body(E)$ in certain positions.
Since barrier events are increasing in $\cM$, compound barrier events 
of $\cM^{*}$ depend only on $Z(\Body(E))$, in an increasing way.

 \end{enumerate}

Condition~\ref{cnd.distr}.\ref{i.distr.indep.hole} says that
for a hole value $E$, the event $\setof{E \in \cH}$
is a decreasing function of $Z(\Body(E))$.
To check this, consider all possible holes of $\cM^{*}$.
We have the following kinds:
 \begin{enumerate}[--]

  \item
Heavy hole values $E$ of $\cM$: all holes of the heavy type of $\cM$
remained holes of $\cM^{*}$, so the event still depends only on
$Z(\Body(E))$ in a decreasing way.

  \item
Hole values $E$ of $\cM^{*}$ of the emerging type: 
such a hole is defined by the property that $\Body(E)$ is a jump.
The jump property is a decreasing function of $Z(\Body(E))$, and therefore
so is the property of being a hole of the emerging type.

  \item
Hole values $E$ of some compound type.
Such a hole appears in $\cM^{*}$ if two holes and a jump
appear in $\cM$ in $\Body(E)$ in certain positions.
Since hole and jump events are decreasing functions of $Z(\Body(E))$
in $\cM$, compound barrier events 
of $\cM^{*}$ depend only on $Z(\Body(E))$, in an increasing way.

 \end{enumerate}

Condition~\ref{cnd.distr}.\ref{i.distr.indep.clean} says first that
for every point $x$ and integer $r$, the events
$\setof{(x,-r) \in \cS}$, $\setof{(x,-r) \in \cS}$  
are decreasing functions of 
$Z((x-r,x))$ and $Z((x,x+r))$ respectively.
The property that $x$ is strongly left $r$-clean in in $\cM^{*}$ 
is defined in terms of strong left-cleanness of $x$ in $\cM$ and the
absence of certain barriers in $[x-r,x]$.
Strong left $r$-cleanness in $\cM$ is a decreasing function of
$Z((x-r,x))$, so is the absence of barriers in $[x-r,x]$.
Therefore strong left $r$-cleanness of $x$
in $\cM^{*}$ is a decreasing function
of $Z((x-r,x))$.

Since both strong and regular left $r$-cleanness in $\cM$ are decreasing
functions of $r$, and the properties stating the
absence of barriers/walls are
decreasing functions of $r$, both strong and regular left $r$-cleanness are
also decreasing functions of $r$ in $\cM^{*}$.
The inequality $\f/3 + \bub < \bub^{*}$, implies that
these functions reach their minimum for $r=\bub^{*}$.
Similar relations hold for right-cleanness.
 \end{proof}

 \begin{lemma}\label{l.cover}
The new mazery $\cM^{*}$ defined by the above construction satisfies
Conditions~\ref{cnd.distr}.\ref{i.distr.inner-clean} 
and~\ref{cnd.distr}.\ref{i.distr.cover}.
 \end{lemma}

 \begin{proof}
Let $E_{1}, E_{2}, \dotsc$ be the sequence of 
maximal external intervals of $\cM$, of
size $\ge \f/3 (> \bub)$.
(We consider $(-\infty,\,0) \sbsq E_{1}$, so that $E_{1}$ is
automatically of size $\ge \f$.)
Let $I_{1},I_{2},\dotsc$ be the intervals betwen them.
By Condition~\ref{cnd.distr}.\ref{i.distr.combinat} of $\cM$, each $I_{j}$
can be covered by a sequence of neighbors $W_{jk}$ in $\cM$.
Every wall of $\cM$ intersects an element of this sequence.
Each pair of these neighbors will be closer than $\f$ to each other.
Indeed, each point of the hop between them belongs either to a 
wall intersecting one of the neighbors, or to a maximal external interval
of size $\le \f/3$, so the distance between the neighbors is at most
$2\bub + \f/3 \le \f$.

Now operation~\ref{i.inform.emerg} above puts some new walls of the
emerging type between the intervals $I_{j}$ or into some of the hops,
all of them disjoint from all existing walls and each other.
If one of these new
walls $W$ comes closer than $\f$ to some interval $I_{j}$ then we add
$W$ to the sequence $W_{jk}$; otherwise, we start a new sequence with $W$.
If as a result of these additions (or originally),
some of these sequences come closer
than $\f$ to each other then we unite them.
By the properties of $\cM$ and the construction of walls of emerging type,
the external intervals between the sequences are hops.
Between the resulting sequences, the distance is $> \f$.
Within these new sequences, every pair of neighbors is closer than $\f$.

Consider one of the above sequences, 
let us call its elements $W_{1},W_{2},\dotsc$.
If it consists of a single light wall $W_{1}$ then it is farther than
$\f$ from all other sequences and operation~\ref{i.inform.rm-light}
removes it.
Since $W_{1}$ is surrounded by maximal intervals, any potential heavy wall
$W$ of $\cM$ intersecting $W_{1}$ is contained in $W_{1}$ and the same
operation removes $W$ as well.

Let us show that the operations of forming compound walls can be used to
create a sequence of consecutive neighbors $W'_{i}$ of $\cM^{*}$
spanning the same interval as $W_{1}, W_{2},\dotsc$.
Assume that walls $W_{i}$ for $i < j$ have been processed already,
and a sequence of neighbors $W'_{i}$ for $i < j'$ has been created
in such a way that 
 \[
   \bigcup_{i<j} W_{i} \sbs \bigcup_{i<j'} W'_{i},
 \]
and $W_{j}$ is not a light wall which is the last in the series.
(This condition is satisfied when $j=1$; indeed, each light wall $W_{i}$
is part of some new wall via one of the compounding operations, since the
ones that are not, would have been removed in
operation~\ref{i.inform.rm-light}.)
We show how to create $W'_{j'}$.

If $W_{j}$ is the last element of the series then it is heavy, and we set
$W'_{j'}=W_{j}$.
Suppose now that $W_{j}$ is not last.

Suppose that it is heavy.
If $W_{j+1}$ is also heavy, or light but not last then $W'_{j'}=W_{j}$.
Else $W'_{j'} = W_{j}+W_{j+1}$, 

Suppose now that $W_{j}$ is light: then it is not last.
If $W_{j+1}$ is last or $W_{j+2}$ is heavy then
$W'_{j'} = W_{j} + W_{j+1}$.
Suppose that $W_{j+2}$ is light.
If it is last then $W'_{j'} = (W_{j} + W_{j+1}) + W_{j+2}$;
otherwise, $W'_{j'} = W_{j} + W_{j+1}$.
 \end{proof}

For condition~\ref{cnd.reachable}.\ref{i.reachable.hop}, in $\cM^{*}$, the
following lemma is needed.

 \begin{lemma}\label{l.new-hop}
Suppose that interval $I$ is a hop of $\cM^{*}$.
Then it is either also a hop of $\cM$ or it contains a sequence
$W_{1},\ldots,W_{n}$ of light walls of $\cM$ separated from each other 
by hops of $\cM$ of size $\ge \f$, and from the ends by hops of $\cM$ of
size  $\ge \f/3$.
 \end{lemma}
 \begin{proof}
If $I$ contains no walls of $\cM$ then it is a hop of $\cM$ and we are done.
Let $U$ be the union $U$ of all walls of $\cM$ in $I$.
The inner cleanness of $I$ in $\cM^{*}$ implies that $U$ is
farther than $\f/3$ from its ends.
Condition~\ref{cnd.distr}.\ref{i.distr.combinat} applied to $U$ implies that
$U$ is spanned by a sequence of neighbor walls $W_{1},W_{2},\dots$ of
$\cM$.
Since $I$ contains no walls of $\cM^{*}$ (and thus no compound walls), 
these neighbor walls are farther than $\f$ from each other.
None of these walls $W_{i}$ is a wall of $\cM^{*}$ therefore each is
contained in an outer clean light wall $W'_{i}$. 
The sequence $W'_{i}$ satisfies our condition: its members are still
separated by hops of size $\ge\f$, for the same reason as $W_{i}$ were.
 \end{proof}

\begin{lemma}\label{l.clean}
The new mazery $\cM^{*}$ defined by the above construction satisfies
Condition~\ref{cnd.distr}.\ref{i.distr.clean}.
\end{lemma}

\begin{proof}
We will use
 \begin{equation}\label{e.bub-f}
  6 \bub < \f,
 \end{equation}
which follows from~\eqref{e.slb-ub}, \eqref{e.g-intro}
and~\eqref{e.f-intro}.
Consider an interval $I$ of size $3 \bub^{*}$ containing no walls of
$\cM^{*}$.
Let $I'$ be the middle third of $I$.
By Condition~\ref{cnd.distr}.\ref{i.distr.inner-clean},
it is contained in a hop of $\cM^{*}$.
Lemma~\ref{l.new-hop} implies that $I'$ is covered by a sequence
$W_{1},\ldots,W_{n}$ of light neighbor 
walls of $\cM$ separated from each other 
by hops of $\cM$ of size $\ge \f$, and surrounded
by hops of $\cM$ of size  $\ge \f/3$.
By~\eqref{e.f-bub-star} we have
$|I'| > 2\f + \bub$, and removing the $W_{i}$ from $I'$ leaves a
subinterval $(a,b) \sbsq I'$ of size at least $\f$.
(If at least two $W_{i}$ intersect $I'$ take the interval between
consecutive ones, otherwise $I'$ is divided into two pieces of total length
at least $2\f$.)
Now $K = (a+\bub+\f/3, b-\bub-\f/3)$ is an interval of length at least
$\f/3-2\bub > 3\bub$ which has distance at least $\f/3$ from any wall.
There will be a clean point in the middle of $K$ which will then be clean
in $\cM^{*}$.
 \end{proof}

Let us look at the reachability conditions~\ref{cnd.reachable}.
Condition~\ref{cnd.reachable}.\ref{i.reachable.wall} will be satisfied for
walls of the emerging type by the (easy) Lemma~\ref{l.emerging-hole}.
For compound walls, it will be shown to be be satisfied
in Subsection~\ref{ss.compound}.

 \begin{lemma}\label{l.if-not-emerg}
Consider the mazery in the middle of the scaling-up operation, after the
construction of all the barriers, walls and holes of the emerging type.
Let $[a, a+\g]$ be an interval contained in a hop of $\cM$ that contains no
emerging walls.
Then for all light barrier types $\ag$ of $\cM$,
some hole of type $\ag'$ is cleanly contained in  $[a+2\bub, a+\g-2\bub]$.  
 \end{lemma}
 \begin{proof}
Suppose that this is not the case.
Then the operation of creating emerging barriers would turn every interval
of the form $[a+\bub+x, a+\g-\bub-y]$ with $0 \le x,y \le \bub$
into an emerging barrier.
Both $a+\bub+x$ and $a+\g-\bub-y$ can be chosen to be clean.
This choice would define an emerging wall.
Since we assumed that we are at the point of the construction when no more
emerging wall can be added, this is not possible.
 \end{proof}

To check condition~\ref{cnd.reachable}.\ref{i.reachable.hop}, we will
proceed as follows. 
Let $u = (x_{0}, x_{1})$ and $v = (y_{0}, y_{1})$ 
be points of $\cV^{*}$ with $\minslope(u, v) \ge \slb^{*}$, such that
for $\d = 0, 1$, the interval $(x_{\d}, y_{\d})$ is a hop in $\cM^{*}$.
We need to prove $u\leadsto v$ in $\cG^{*}$.
Since the intervals $(x_{\d}, y_{\d})$ are hops in $\cM^{*}$,
Lemmas~\ref{l.if-not-emerg} and \ref{l.new-hop} show that they satisfy the
conditions of Lemma~\ref{l.grate}.

\section{Emerging and compound types}\label{s.bounds}

In this section, we give bounds on the probabilities of emerging and
compound barriers and holes, as much as this is possible without indicating
the dependence on $k$.

\subsection{General bounds}\label{ss.bounds}

We start with some estimates that are needed for both kinds of types.
Let $\pub$ be some upper bound on the sum of all barrier probabilities:
 \begin{equation}\label{e.pub}
 \begin{split}
            \pub &\ge \sum_{\r \ge \R} \p(\r).
 \end{split}
 \end{equation}
Let $\ag$ be a barrier type, $a \le b < c$ and $(b-a), (c-b)  \le 6 \bub$,
and let $E(a, b, c, \ag)$ be defined as in
Condition~\ref{cnd.distr}.\ref{i.distr.hole-lb}.
We extend the bound~\eqref{e.hole-lb} of this condition in several ways.
These extension lemmas rely explicitly on
Condition~\ref{cnd.distr}.\ref{i.distr.hole-lb}.
For the following lemma, remember the definition of $\ncln$
before~\eqref{e.ncln-ub}.

 \begin{lemma}\label{l.hole-lb-clean}
Let $F(a,b,c,\ag)$ be the event that $E(a,b,c,\ag)$ occurs and also
the end of the hole is strongly right-clean.
 \begin{equation}\label{e.hole-lb-clean}
   \Prob(F(a,b,c,\ag)) \ge (1 - \ncln) \Prob(E(a,b,c,\ag)).
 \end{equation}
 \end{lemma}
 \begin{proof}
For $b \le x \le c+\bub$, let $E_{x}$ be the event that
$E(a,b,c,\ag)$ is realized by a hole ending at $x$ but is not realized
by any hole ending at any $y < x$.
Let $F_{x}$ be the event that $x$ is strongly right-clean.
Then $E(a,b,c,\ag) = \bigcup_{x}E_{x}$,
$F(a,b,c,\ag) \sps \bigcup_{x} (E_{x} \cap F_{x})$,
the events $E_{x}$ and $F_{x}$ are independent for each $x$, 
and the events $E_{x}$ are mutually disjoint. 
Hence
 \begin{equation*}
  \begin{split}
   \Prob(F(a,b,c,\ag)) &\ge \sum_{x} \Prob(E_{x}) \Prob(F_{x})
                   \ge (1-\ncln) \sum_{x} \Prob(E_{x}) 
\\                 &= (1-\ncln) \Prob(E(a,b,c,\ag)).
  \end{split}
 \end{equation*}
 \end{proof}

We extend the hole lower bound~\eqref{e.hole-lb} here to cases when
$b - a > 6 \bub$, though with a different constant.

 \begin{lemma}\label{l.hole-lb-1}
Let $a \le b < c$ with $c-b  \le 6 \bub$: then we have
 \[
   \Prob(E(a,b,c,\ag)) \ge 
      (1 - \ncln - (c-a)\pub) (c-b)^{\hxp} \h(\ag).
 \]
 \end{lemma}
 \begin{proof}
If $b - a \le 6 \bub$ then we can apply the hole lower
bound~\eqref{e.hole-lb}; suppose therefore that this does not hold.
Let $a' = b-\bub$.
Then~\eqref{e.hole-lb} is applicable to $(a',b,c)$, and we get
 \[
   \Prob(E(a',b,c,\ag)) \ge (c-b)^{\hxp}\h(\ag).
 \]
Consider the event $C$ that $a$ is strongly right-clean
and the interval $(a,c)$ contains no barriers.
Then $C \cap E(a',b,c,\ag) \sbs E(a,b,c,\ag)$.
Event $C$ is decreasing with $\Prob(C) \ge 1 -\ncln - (c-a)\pub$.
By the FKG inequality, we have 
$\Prob(E(a,b,c,\ag)) \ge \Prob(C)\Prob(E(a',b,c,\ag))$. 
 \end{proof}

Finally, we extend the hole lower bound~\eqref{e.hole-lb} to cases when 
$c - b > 6 \bub$.

 \begin{lemma}\label{l.hole-lb-2}\
Let $a \le b < c$.
  \begin{enumerate}[\upshape 1.]

   \item\label{i.hole-lb-2.main}
We have
 \begin{equation}
   \Prob(E(a, b, c, \ag)) \ge 
    (1 - \ncln - (c-a)\pub) (0.6 \et (0.1 (c-b)^{\hxp} \h(\ag))).
 \end{equation}

   \item\label{i.hole-lb-2.higher}
If the event $E^{*}(a, b, c, \ag)$ is defined like $E(a, b, c, \ag)$ except
that the jump in question must be a jump of $\cM^{*}$ then,
assuming $\pub^{*} \le \pub$, we have
 \begin{equation}\label{e.hole-lb-2}
   \Prob(E^{*}(a, b, c, \ag)) \ge 
 (1 - \ncln^{*} - 2(c-a)\pub) (0.6 \et (0.1 (c-b)^{\hxp} \h(\ag))).
 \end{equation}

  \end{enumerate}
 \end{lemma}
 \begin{proof}
We will use the following inequality, which can be checked by direct
calculation.
Let $v = 1 - 1/e = 0.632\dots$, then for $x > 0$ we have
 \begin{equation}\label{e.exp-lb}
   1 - e^{-x} \ge v \et v x.
 \end{equation}
In view of Lemma~\ref{l.hole-lb-1}, for
the first statement of the lemma, we
only need to consider the case $c - b > 6 \bub$.

Let $n = \Flo{\frac{c - b}{3\bub}}-1$, then we have
 \begin{equation}\label{e.n-lb}
  n \bub \ge (c-b)/9.
 \end{equation}
Let
 \[
 \begin{split}
     a_{i} &= b + 3 i \bub,
\quad  E_{i} = E(a_{i}, a_{i} + \bub, a_{i} + 2\bub, \ag),
\quad \txt{ for } i=1,\dots,n,
\\     E  &= \bigcup_{i} E_{i},
\quad      \s = \bub^{\hxp}\h(\ag).
 \end{split}
 \]
Inequality~\eqref{e.hole-lb} is applicable to
$E_{i}$ and also $s \le 0.6$ by~\eqref{e.0.6}.
We have $\Prob(E_{i}) \ge \s$, hence
$\Prob(\neg E_{i}) \le 1 - \s \le e^{-\s}$.
The events $E_{i}$ are independent, so 
 \begin{equation}\label{e.union-E}
  \Prob(E) = 1 - \prod_{i} \Prob(\neg E_{i}) 
                           \ge 1 - e^{- n \s} \ge 0.6 \et (0.6 n \s),
 \end{equation}
where in the last step we used~\eqref{e.exp-lb}.
We have
 \begin{equation}\label{e.ns}
    n \s = n \bub^{\hxp} \h(\ag).
 \end{equation}
Now, by~\eqref{e.n-lb}, we have 
$n \bub^{\hxp} \ge (\bub n)^{\hxp} \ge 9^{-\hxp}(c-b)^{\hxp}$.
Substituting into~\eqref{e.union-E}, \eqref{e.ns}:
 \[
 \begin{split}
   \Prob(E) 
    &\ge 0.6 \et (0.6 \cdot 9^{-\hxp} (c-b)^{\hxp}\h(\ag))
    \ge 0.6 \et (0.1 (c-b)^{\hxp} \h(\ag)),
 \end{split}
 \]
where we used $\hxp < \log_{9} 6$, which follows from~\eqref{e.hxp}.
The event $E$ implies that a hole of type $\ag'$
starts in $[b,c-1]$ and that the left end of the hole is strongly left-clean.
Let $C$ be the event that $a$ is strongly right-clean and
that there is no barrier in $(a,c)$.
If also $C$ holds then there is a jump from $a$ to our hole, which we
need.
The event $C$ is decreasing, so the FKG inequality implies that
$\Prob(C) \ge 1 - \ncln - (c-a)\pub$
can be multiplied with $\Prob(E)$ for a lower bound.
For the second statement, the requirements must be added that
there are no barriers of $\cM^{*}$ in $(a,c)$, and that $a$ must
be strongly right-clean in $\cM^{*}$.
We can replace $\ncln$ with the larger $\ncln^{*}$,
and we can replace $\pub$ with $2\pub \ge \pub + \pub^{*}$.
 \end{proof}

 \begin{remark}
It may seem that the proof of
Lemmas~\ref{l.hole-lb-1} and \ref{l.hole-lb-2} would go through
even with $(c-b)\h(\ag)$ in place of $(c-b)^{\hxp}\h(\ag)$.
However, we used inequality~\eqref{e.0.6} (the smallness
needed for the approximation by $e^{x}$).
Even if we assume that this inequality holds for $\cM$ without $\hxp$, we
can prove it for $\cM^{*}$ only with $\hxp$.
 \end{remark}

\subsection{The emerging type}\label{ss.emerg}

Recall the definition of an emerging type in
part~\ref{i.inform.emerg} of the scale-up algorithm in
Section~\ref{s.plan}.
Recall also that dimension 0 is the direction of the $X$ sequence
and dimension $1$ the direction of the $Y$ sequence, thus $\cW_{0}$ is the
set of walls defined by $X$.

 \begin{lemma}\label{l.emerging-hole}
If a wall of the emerging type $\ang{\g}$ and size $w_{1}$
begins at $x$ in $\cW_{0}$ and a hole of type
$\ang{\g}'$ and size $w_{2}$ begins at $y$ in $\cH_{1}$ then in graph $\cG$
there is a path of slope $\le 1$ from $(x,y)$ to $(x + w_{1}, y + w_{2})$.
 \end{lemma}
 \begin{proof}
This follows directly from the reachability
condition~\ref{cnd.reachable}.\ref{i.reachable.hop},
if we observe that $w_2 \le w_1 \le 2w_2$, so that the the minslope is 
$\ge 1/2$.
 \end{proof}

 \begin{lemma}\label{l.emerging}
Let $n = \Flo{\frac{\g - 5\bub}{3 \bub}}$.
For any point $x$, the expression
 \begin{equation}\label{e.emerging-wall}
   2\bub |\Btypes| e^{- (1 - \ncln) n \h(\R^{*})}
 \end{equation}
is an upper bound on the sum, over all $w$, of the
probabilities  that an emerging barrier of type $\ang{\g}$ (with rank
$\R^{*}$) and size $w$ starts at $x$. 
 \end{lemma}
 \begin{proof}
Recall the definition of an emerging barrier.
Suppose that there are $a \le a' < b' \le b$ with
$a'-a \le \bub$, $b'-a' = \g-4\bub$, $b-b' \le \bub$ and a
light barrier type $\ag$ such that 
no hole of type $\ag'$ is cleanly contained in $[a',b']$. 
Then $((a,b),\ang{\g})$ is a barrier of $\cM^{*}$.
This definition implies that if $(x,y)$ is an emerging barrier then there
is a light barrier type $\ag$ such that no hole of type $\ag'$ is cleanly
contained in $[x+\bub,x+\g-5\bub]$.

Let us fix $\ag$, and for any
$i = 0,\dots,n-1$, let $\cA(i, \ag)$ be the event that
no strongly outer-clean hole of type $\ag'$ starts at $x + (3 i + 1) \bub$.
For each fixed $\ag$, these events are independent.
Indeed, according to Condition~\ref{cnd.distr}.\ref{i.distr.indep},
event $\cA(i, \ag)$ only depends on $Z((x + 3 i\bub), (x + 3(i+1)\bub))$.

Recall the hole lower bound:
Condition~\ref{cnd.distr}.\ref{i.distr.hole-lb}.
It says the following.
For $a \le b < c$ and $b-a, c-b  \le 6 \bub$,
and let $E(a, b, c, \ag)$ be the event
that there is a $d \in [b, c-1]$ such that $[a,d]$ is a jump
and a hole of type $\ag'$ starts at $d$.
Then 
 \[
   \Prob(E(a,b,c,\ag)) \ge (c-b)^{\hxp} \h(\ag).
 \]
With $a=x+3i\bub$, $b=x+(3i+1)\bub$, $c=x+(3i+1)\bub+1$, this implies
that with probability at least $\h(\ag) \ge \h(\R^{*})$, 
a strongly left-clean hole starts at $x+(3i+1)\bub$.
Lemma~\ref{l.hole-lb-clean} implies that with probability
at least $\h(\R^{*})(1-\ncln)$, this hole is also strongly right-clean, and
so it is strongly outer-clean.
Hence, each of the independent events $\cA(i, \ag)$ has 
a probability upper bound $1 - (1 - \ncln) \h(\R^{*})$.
The probability that all the events $\cA(i,\ag)$ hold is bounded by 
$(1 - (1 - \ncln) \h(\R^{*}))^{n} \le e^{-n (1 - \ncln) \h(\R^{*})}$.
The probability that one of these events holds for some $\ag$ is
at most $|\Btypes|$ times larger.
The sizes of emerging barriers vary in the range $[\g-4\bub,\g-2\bub]$,
hence the factor $2\bub$.
 \end{proof}

\subsection{Compound types}\label{ss.compound}

 \begin{lemma}\label{l.compound-hole}
For $i=1,2$, let $W_{i}$ be two neighboring walls of $\cW_{0}$
starting at points $x_{1} < x_{2}$, with sizes $w_{i}$ with distance $d$.
Let $H_{i}$, for $i=1,2$ be two disjoint holes of $\cH_{1}$ with sizes
$h_{i}$, with starting points $y_{i}$ where $H_{i}$ fits $W_{i}$.
Suppose that the interval between these holes is a hop.
Let $e$ be the distance between $H_{i}$.
 \begin{alphenum}

  \item 
The reachability relation
$(x_{1}, y_{1})\leadsto (x_{2} + w_{2}, y_{2} + h_{2})$
is implied by the inequality
  \begin{equation}\label{e.e-cond}
   \slb d \le e \le d.
 \end{equation}
  
  \item\label{i.compound-fitting}
Let $\dlb < \dub$, $\elb < \eub$, be integers satisfying
 \begin{equation}\label{e.compound-fitting}
   \slb (\dub - 1) \le \elb \le \eub - 1 \le \dlb.
 \end{equation}
  If $\dlb \le d < \dub$ and $\elb \le e < \eub$ then 
inequality~\eqref{e.e-cond} is satisfied.

 \end{alphenum}
 \end{lemma}

 \begin{proof}
Assume~\eqref{e.e-cond}.
As $W_{1}$ is passed by $H_{1}$, there is a path from
$u_{1} = (x_{1},y_{1})$ to 
 \[
  u_{2} = (x_{1} + w_{1},\; y_{1} + h_{1}).
 \]
Due to~\eqref{e.e-cond} there is a path from $u_{2}$ to  
$u_{3} = (x_{2}, y_{2})$.
As $W_{2}$ is passed by $H_{2}$, there is a path from $u_{3}$ to 
$(x_{2} + w_{2},\; y_{2} + h_{2})$.
The total slope of the combination of these three paths is 
clearly at most 1.
(For reachability, the total slope condition is not important; but, 
a compound hole will need to satisfy
Condition~\ref{cnd.types}.\ref{i.types.fitting}.)
The proof of statement~\eqref{i.compound-fitting} is immediate.
 \end{proof}

Recall the definition of compound (barrier and wall) types given in
part~\ref{i.inform.compound} of the scale-up algorithm of
Section~\ref{s.plan}.
In~\eqref{e.jump-types}, we defined a sequence of integers $d_{i}$.
For any pair $\ag_{1},\ag_{2}$ of barrier types where $\ag_{1}$ is light,
and any $0 \le i$ with $d_{i} \le \f$, there is a new type
$\bg = \ang{\ag_{1},\ag_{2},i}$ in $\Btypes^{*}$.
A barrier of type $\bg$ occurs in $\cM^{*}$ wherever disjoint
barriers $W_{1}, W_{2}$ of types $\ag_{1}, \ag_{2}$ occur (in this order)
at a distance $d \in [d_{i},d_{i+1}-1]$.
Let the barriers $W_{i}$ have starting points $x_{1} < x_{2}$, and sizes
$w_{i}$.
The body of the new barrier is $(x_{1}, x_{2} + w_{2})$.
The barrier becomes a wall if the component barriers are walls and the
interval between them is a hop.

For compound hole types, we used the sequence $e_{i}$ defined
in~\eqref{e.hole-jump-types}.
A hole of type $\bg'$ occurs in $\cM^{*}$
if holes $H_{1}, H_{2}$ of types $\ag'_{1},\ag'_{2}$ occur,
connected by a jump of $\cM$ of size $d \in [e_{i}, e_{i+1}-1]$.

 \begin{lemma}\label{l.compound-hole-fitting}
For each $i$ if we set 
$\dlb=d_{i}$, $\dub=d_{i+1}$, $\elb=e_{i}$, $\eub=e_{i+1}$ then
inequality~\eqref{e.compound-fitting} is satisfied.
Therefore for each compound type $\bg$, if $W$ is a vertical compound wall
of type $\bg$ with body $(x_{1},x_{2})$ and and $H$ is a horizontal
wall of type $\bg'$ with body $(y_{1},y_{2})$ then 
$(x_{1},y_{1})\leadsto(x_{2},y_{2})$: the hole ``passes''
through the wall.
 \end{lemma}
 \begin{proof}
Assume $i < 16$.
then $d_{i}=i=e_{i}$, and $d_{i+1}=e_{i+1}=i+1$.
Inequalities~\eqref{e.compound-fitting} turn into the true inequalities
$\slb i \le i \le i \le i$.
Assume $i=16$, then $d_{i}=e_{i}=16$, and $d_{i+1}=19$, $e_{i+1}=17$.
Inequalities~\eqref{e.compound-fitting} will turn into the true inequalities
$\slb \cdot 18 \le 16 \le 16 \le 16$.
Assume $i=17$, then $d_{i}= 19$, $d_{i+1}=22$, $e_{i}=17$, $e_{i+1}=19$.
Inequalities~\eqref{e.compound-fitting} will turn into the true inequalities
$\slb\cdot 21 \le 17 \le 18 \le 19$.
Assume $i>17$, then $d_{i}=\flo{\lg^{i}}$, $e_{i}=\flo{\lg^{i-1}}$.
Given $\slb \le \half$ and $\lg = 2^{1/4}$, what we need to check 
from~\eqref{e.compound-fitting} is
 \[
 \frac{1}{2}(\flo{\lg^{i+1}} - 1) \le \flo{\lg^{i-1}} \le \flo{\lg^{i}}-1,
 \]
which is true for $i > 17$.
 \end{proof}

 \begin{lemma}\label{l.compound-wall-ub}
For any $\r_{1},\r_{2}$, the sum, over all $w$, of the probabilities
for the occurrence of a compound
barrier of type $\ang{\ag_{1},\ag_{2},i}$ with $\Rank(\ag_{j})=\r_{j}$
and width $w$ at a given point $x_{1}$ is bounded above by
 \begin{equation}\label{e.compound-wall-ub}
   (d_{i+1} - d_{i}) \p(\r_{1}) \p(\r_{2}).
 \end{equation}
 \end{lemma}
 \begin{proof}
For fixed $\r_{1},\r_{2},x_{1},d$, let $B(d,w)$ be the event that a compound
barrier of any type $\ang{\ag_{1},\ag_{2},i}$ with
$\Rank(\ag_{j})=\r_{j}$, distance $d$ between the component barriers,
and size $w$ appears at $x_{1}$.
For any $w$, let $A(x,\r,w)$ be the event that a barrier of rank
$\r$ and size $w$ starts at $x$.
We can write
 \[
   B(d,w) = \bigcup_{w_{1}+d+w_{2}=w} 
          A(x_{1},\r_{1},w_{1}) \cap A(x_{1}+w_{1}+d,\r_{2},w_{2}).
 \]
where events $A(x_{1},\r_{1},w_{1})$, $A(x_{1}+w_{1}+d,\,\r_{2},w_{2})$ are
independent.
By~\eqref{e.p(,)}:
 \[
   \Prob(B(d,w)) \le \sum_{w_{1}+d+w_{2}=w} 
          \p(r_{1},w_{1}) \p(\r_{2},w_{2}).
 \]
Hence by~\eqref{e.p()}:
 \[  
 \sum_{w} \Prob(B(d,w)) \le \sum_{w_{1}} \p(\r_{1},w_{1})
                  \sum_{w_{2}} \p(\r_{2},w_{2}) \le \p(\r_{1})\p(\r_{2}).
 \]
 \end{proof}

 \begin{lemma}\label{l.compound-hole-lb}
Consider the compound hole type $\bg'$ where
$\bg = \ang{\ag_{1}, \ag_{2}, i}$.
For $a \le b < c$, let $E_{2}(a, b, c, \bg)$ be the event
that there is a $d \in [b, c-1]$ such that $[a,d]$ is a
jump of $\cM^{*}$, and a compound hole of type $\bg$ starts at $d$.
Assume
 \begin{equation}\label{e.s-ub}
                                 c - a \le 12 \bub^{*},
\quad    (\bub^{*})^{\hxp} \h(\ag_{i}) \le 0.5,
\quad                         \pub^{*} \le \pub.
 \end{equation}
Then we have
 \begin{equation}\label{e.compound-hole-lb}
 \begin{split}
   \Prob(E_{2}(a,b,c,\bg)) 
\ge 0.01 (1 - 2\ncln^{*} - 25 \bub^{*} \pub) 
  (c-b)^{\hxp} (e_{i+1} - e_{i})^{\hxp} \h(\ag_{1}) \h(\ag_{2}).
 \end{split}
 \end{equation}
 \end{lemma}
 \begin{proof}
Let $\elb=e_{i}$, $\eub = e_{i+1}$.
For each $x \in [b, c + \bub - 1]$, 
let $A_{x}$ be the event that there is a $d_{1} \in [b, c-1]$ such that
$[a, d_{1}]$ is a jump of $\cM^{*}$, and a hole of type $\ag'_{1}$ starts at
$d_{1}$ and ends at $x$, and that $x$ is the smallest possible number with
this property.
Let $B_{x}$ be the event that there is a $d_{2} \in [x + \elb, x+\eub)$
such that 
$[x, d_{2}]$ is a jump of $\cM$, and a hole of type $\ag'_{2}$ starts at
$d_{2}$.
Then $E_{2}(a,b,c,\bg) \sps \bigcup_{x} (A_{x} \cap B_{x})$,
and for each $x$, the events $A_{x},B_{x}$ are independent.
We have, using the notation of Lemma~\ref{l.hole-lb-2}:
 \[
 \begin{split}
   \sum_{x} \Prob(A_{x}) &= \Prob(E^{*}(a, b, c, \ag_{1})) 
\\   &\ge
 (1 - \ncln^{*} - 2(c-a)\pub) (0.6 \et (0.1 (c-b)^{\hxp} \h(\ag_{1}))).
 \end{split}
 \]
Further, using the same lemma:
 \[
   \Prob(B_{x}) = \Prob(E(x, x + \elb, x + \eub, \ag_{2}))
  \ge (1 - \ncln - \eub \pub) 
      (0.6 \et (0.1(\eub - \elb)^{\hxp} \h(\ag_{2}))).
 \] 
By the assumptions~\eqref{e.s-ub}:
$0.1 (c-b)^{\hxp} \h(\ag_{1}) 
\le 0.1(12\bub^{*})^{\hxp} \h(\ag_{1}) \le 0.6$,
hence the operation $0.6 \et$ can be deleted.
The same reasoning applies to the second application of $0.6\et$.
Combining these, using $\ncln^{*} > \ncln$, $\pub^{*} < \pub$:
 \[
  \begin{split}
   \Prob(E_{2}(a, b, c, \bg)) &\ge \sum_{x} \Prob(A_{x})\Prob(B_{x})
\\ &\ge (1 - 2\ncln^{*} - (2(c-a) + \eub)\pub) 0.01
   (c-b)^{\hxp} (\eub - \elb)^{\hxp} \h(\ag_{1}) \h(\ag_{2}).
  \end{split}
 \]
 \end{proof}

\section{The scale-up functions}\label{s.params}

\subsection{Parameters}

Lemma~\ref{l.main} says, with $\p = \Prob\setof{Z(i)=1}$,
that there is a $\p_{0}$ such that if $\p < \p_{0}$ then the
sequence $\cM^{k}$ can be constructed in such a way that~\eqref{e.main}
holds. 
Our construction has several parameters.
If we computed $\p_{0}$ explicitly then all these parameters could be
turned into constants: but this is unrewarding work and it would only make
the relationships between the parameters less intelligible.
We prefer to name all these parameters, point out the necessary
inequalities among them, and finally show that if $\p$ is sufficiently
small then all these inequalities can be satisfied simultaneously.

Recall that the slope lower bound $\slb$ must satisfy
 \begin{equation}\label{e.slb-lb}
   \slb < \slb_{0} = 1/2.
 \end{equation}
The parameter $\lg > 0$ was introduced in~\eqref{e.lambda-def}.
We will use a parameter $\R$ (not necessarily integer) 
for a lower bound on all ranks in the mazery.

It can be seen from the definition of compound ranks
in~\eqref{e.compound-rank} and from
Lemma~\ref{l.compound-wall-ub} that the probability bound $\p(\r)$ 
of a barrier type should be approximately $\lg^{-\r}$.
We introduce an upper bound that is a little smaller:
 \begin{equation}\label{e.wall-prob}
  \p(\r) = \aux_{2}\r^{-\aux_{1}} \lg^{-\r}
 \end{equation}
where $\aux_{1}$ has been defined in~\eqref{e.Rlb.aux1} above, and 
 \begin{equation}\label{e.aux2}
 0 < \aux_{2} = 0.15 < 1.
 \end{equation}
The inequality requiring this definition is~\eqref{e.Rlb.aux2.pub} below.
The term $\aux_{2}\r^{-\aux_{1}}$, just like the factor
in the function $\h(\r)$ defined in the hole lower bound~\eqref{e.hole-lb},
serves for absorbing some lower-order factors that arise in estimates
like~\eqref{e.compound-hole-lb}.
We define $\aux_{3}$ and then $\aux_{0}$ by
 \begin{equation}\label{e.Rlb.aux3}
  \aux_{0}\aux_{2}^{\hxp} = \aux_{3} = 436,
 \end{equation}
this gives
 \[
  \h(\r) = \aux_{3}\lg^{-\r\hxp}.
 \]
The value of $\aux_{3}$ is required by the 
inequality~\eqref{e.Rlb.aux3-lb} below.
This defines $\aux_{0}$ implicitly
using the values of $\hxp$ from~\eqref{e.hxp} and $\aux_{2}$
from~\eqref{e.aux2}.

By these definitions, we can give a concrete value to the 
upper bound $\pub$ introduced in~\eqref{e.pub}:
 \begin{equation}\label{e.pub-def}
   \pub = \lg^{-\R},
\quad \T = 1/\pub = \lg^{\R}
 \end{equation}
where we have also introduced its inverse, $\T$.
 
 \begin{lemma}\label{l.pub}
The above definition of $\pub$ satisfies~\eqref{e.pub}.
 \end{lemma}
  \begin{proof}
We have
 \begin{equation}\label{e.Rlb.aux2.pub}
  \sum_{\r \ge \R} \p(\r) <
  \aux_{2} \sum_{\r \ge \R} \lg^{-\r}
  = \lg^{-\R}\frac{\aux_{2}}{1 - 1/\lg} < \lg^{-\R}
\quad\txt{ if }\quad
   \aux_{2} < 1 - 1/\lg,
 \end{equation}
which is satisfied by the choice $\aux_{2}$ in~\eqref{e.aux2}.
  \end{proof}

Several other parameters of $\cM$ and the scale-up are expressed
conveniently in terms of $\T$:
 \begin{equation}\label{e.bubxp-etc-def}
                       \bub      = \T^{\bubxp},
\quad                  \f        = \T^{\fxp},
\quad                  \g        = \T^{\gxp},
\quad 0 < \bubxp < \gxp < \fxp < 1.
 \end{equation}
To obtain the new rank lower bound, we multiply $\R$ by a constant:
 \begin{equation}\label{e.txp}
 \begin{split}
           \R = \R_{k}  = \R_{1} \txp^{k-1},
\quad       \R_{k+1} = \R^{*} = \R \txp,
\quad                1 < \txp, \R_{1}.
 \end{split}
 \end{equation}
It is convenient to introduce $\R_{0}=\R_{1}/\txp$, to be able to write
$\R_{k}=\R_{0}\txp^{k}$.
Setting $\R_{0}$ large enough will enable to us to satisfy any inequality
of the form $c\T^{-x} < 1$ for each mazery in the sequence as long as the
constants $c$ and $x$ are strictly positive.
We will collect the required bounds on $\R_{0}$ as we go along.
A few more consequences of these definitions:
 \[
   \T^{*} = \lg^{\R^{*}} = \lg^{\R\txp} = \T^{\txp},
\quad \bub^{*} = \bub^{\txp},
\quad \log_{\lg} \f = \R\fxp.
 \]
A bound on $\txp$ has been indicated in the requirement~\eqref{e.R-cond}
which will be satisfied if
 \begin{equation}\label{e.txp-ub}
   \txp \le 2 - \fxp.
 \end{equation}
Let us make sure that Condition~\ref{cnd.types}.\ref{i.types.bub}
is satisfied by $\cM^{*}$.
Barriers of the emerging type
have size at most $\g-2\bub$, and at the time of their creation, they
are the largest existing ones.
We get the largest new barriers when the compound operation combines these
with light barriers on both sides, leaving the largest gap possible,
so the largest new barrier size is 
$\g - 2\bub + 2(\f + \bub) \le \g+2\f \le 3\f$, where we 
used~\eqref{e.f-intro}, \eqref{e.g-intro}.
Hence any value larger than $3\f$ can be chosen as 
$\bub^{*} = \bub^{\txp}$. 
With $\R_{0}$ large enough, we always get this if
 \begin{equation}\label{e.bubxp-lb}
  \fxp < \bubxp \txp < 1
 \end{equation}
(where the second inequality will also be needed).\footnote{More exactly, we need
$\R \ge 64$.}
We can satisfy~\eqref{e.g-intro} similarly, if
 \begin{equation}\label{e.Rlb.g-intro.final}
  \R \ge 130.
 \end{equation}

The exponent $\hxp$
has been part of the definition of a mazery: it is the power by which,
roughly, hole probabilities are larger than barrier probabilities. 
We require
 \begin{equation}\label{e.hxp-ub}
   0 < \hxp < (\gxp - \bubxp)/\txp.
 \end{equation}

 \begin{lemma}\label{l.params-choice-1}
 The exponents $\bubxp,\fxp,\gxp,\txp,\hxp$ can be chosen to
satisfy the inequalities~\eqref{e.bubxp-etc-def}, \eqref{e.txp},
\eqref{e.txp-ub},  \eqref{e.bubxp-lb}, \eqref{e.hxp-ub}.
 \end{lemma}
 \begin{proof}
We can choose $\hxp$ last, to satisfy~\eqref{e.hxp-ub}, so consider just
the other inequalities.
Choose $\txp = 2 - \fxp$ to satisfy~\eqref{e.txp-ub}; 
then~\eqref{e.bubxp-etc-def} and \eqref{e.bubxp-lb} will be satisfied if 
$\bubxp < \gxp < \fxp < \bubxp (2 - \fxp) < 1$. 
This is achieved by
 \begin{equation}\label{e.Rlb.exps}
 \bubxp = 0.4,\quad \gxp = 0.45,\quad \fxp = 0.5,\quad \hxp = 0.03,
\txt{ hence } \txp = 1.5.
 \end{equation}
 \end{proof}
Let us fix now all these exponents.
In order to satisfy all our requirements also for small $k$,
we fix $\aux_{1}$ next; then we fix $\aux_{0}$ and finally $\R_{0}$.
Each of these last three parameters just has to be chosen sufficiently
large as a function of the previous ones.

We need upper bounds on the largest ranks, and on the number of types.

  \begin{lemma}\label{l.walltypes-ub}\
 \begin{enumerate}[\upshape 1.]

   \item The quantity $\R\frac{2 \txp}{\txp - 1}$
is an upper bound on all existing ranks in a mazery.
Hence every rank exists in $\cM^{k}$ for at most
$\log_{\txp}\frac{2 \txp}{\txp -1}$ values of $k$.

    \item We have, denoting for the moment $c = \log_{\txp} 3$:
 \begin{equation}\label{e.walltypes-ub}
   |\Btypes| < \frac{\R_{0}^{(\R/\R_{0})^{c}}}{\txp\R}
 = \frac{\R_{0}^{3^{k}}}{\txp\R}
 \end{equation}

 \end{enumerate}
  \end{lemma}
 \begin{proof}
For the moment, let us denote the largest existing rank by $\Rmax$.
Emerging types got a rank equal to $\R^{*}$, 
and the largest rank produced by the compound operation is at most
$\Rmax + 2\R^{*}$ (since the compound operation is applied twice),
hence $\Rmax_{k+1} \le \Rmax_{k} + 2\R_{k+1}$.
Since also $\Rmax_{1} \le 2\R_{1}$ (since there is only one rank in
$\cM^{1}$), we have for $k \ge 1$:
 \begin{equation}
 \Rmax_{k} \le 2 \sum_{i=1}^{k} \R_{i} 
                = 2 \R_{0}\txp\frac{\txp^{k}-1}{\txp - 1}
                \le \R_{k}\frac{2 \txp}{\txp - 1}.
 \end{equation}

Now for the number of types.
There is only one emerging type.
The operation of forming compound types multiplies the number of
types at most by the number of values $i$ in the
definition~\eqref{e.jump-types}: this is $\le \f$ for $\f < 17$
and $\le \log_{\lg} \f = \R\fxp$ otherwise.
For the moment, let $N$ denote the number of barrier types.
The operation of forming compound types once results in 
multiplying $N$ by at most $N\R\fxp$ and adding the
result to $N$.
We have to repeat this operation twice, and use $\fxp=0.5$:
 \[
   N^{*} \le 1 + N(1 + N\R\fxp)^{2} = 
   \R^{2}N^{3}(\R^{-2}N^{-3} + \R^{-2}N^{-2} + \R^{-1}N^{-1}+ 0.25)
   <  \R^{2}N^{3}
 \]
if $\R \ge 3$.
This recursive inequality leads to the estimate~\eqref{e.walltypes-ub}.
This is straightforward with the recursion $N^{*} \le N^{3}$ since what we
are proving is $N_{k} \le \R_{0}^{3^{k}-1}/\txp^{k+1}$.
The divisor $\R$ in~\eqref{e.walltypes-ub} absorbs the effect of
the factor $\R^{2}$ in the recursion.
 \end{proof}

 \subsection{Probability bounds after scale-up}

The structures $\cM^{k}$ are now defined but we have not proved yet that
they are mazeries, since not all inequalities required in the definition of
mazeries have been verified yet.

 \begin{lemma}\label{l.0.6}
If $\R_{0}$ is sufficiently large then for each $k$, for the structure
$\cM^{k}$, for any barrier type $\ag$, inequality~\eqref{e.0.6} holds; 
we also have
 \begin{equation}\label{e.sum-bub-k+1}
  \sum_{k} \bub_{k+1} \pub_{k} < 0.5.
 \end{equation}
 \end{lemma}

 \begin{proof}
Let us prove~\eqref{e.0.6}, which says $\bub^{\hxp} \h(\R) < 0.6$.
We have
$\bub^{\hxp}\h(\R) = \T^{\bubxp\hxp}\aux_{3}\T^{-\hxp}
    = \aux_{3} \T^{-\hxp(1 - \bubxp)}
$
which is smaller than 0.6 if $\R_{0}$ is sufficiently large.\footnote{More 
exactly, we need 
$\R \ge 2113$.}
For inequality~\eqref{e.sum-bub-k+1}, note that
 \[
  \sum_{k}\bub_{k+1} \pub_{k} = 
   \sum_{k} \lg^{-\R_{0}\txp^{k}(1 - \bubxp\txp)}
 \]
which because of~\eqref{e.bubxp-lb}, 
is clearly less than 0.5 if $\R_{0}$ is large.\footnote{More exactly,
we need at most $\R_{0} \ge 31$.}
 \end{proof}

Recall the constant $\slopeincr$ defined in~\eqref{e.slopeincr}.
Note that for $\R_{0}$ large enough, the relations
 \begin{align}
\label{e.strong-ncln-ub}
   \bub^{*}\pub &< 0.5(0.25 - \ncln),
\\\label{e.strong-f-intro}
       \slopeincr \g/\f &< 0.5(\slb_{0} - \slb).
 \end{align}
hold for $\cM=\cM^{1}$.\footnote{More exactly, 
since $\ncln_{1}=0$, for~\eqref{e.strong-ncln-ub} we need
$\R \ge 31$.}
Further, since $\slb_{1}=0$, for~\eqref{e.strong-f-intro} we need
 \begin{equation}\label{e.Rlb.g-per-f-start}
 \begin{split}
    \slb_{0}/2 &> \slopeincr \g / \f = \slopeincr \T^{-(\fxp - \gxp)},
 \end{split}
 \end{equation}
satisfied if $\R \ge 1517$.
The following lemma establishes
Condition~\ref{cnd.distr}.\ref{i.distr.bounds.ncln-ub} and
inequality~\eqref{e.f-intro} for all $k$.

\begin{lemma}\label{l.ncln-ub}\
Suppose that the structure $\cM=\cM^{k}$ is a mazery and it
satisfies~\eqref{e.strong-ncln-ub} and~\eqref{e.strong-f-intro}.
Then $\cM^{*}=\cM^{k+1}$ also satisfies these inequalities.
 \end{lemma}

 \begin{proof}
The probability that a point $a$ is strongly clean in $\cM$ but not
in $\cM^{*}$
is clearly upperbounded by $\bub^{*} \pub$, which upperbounds the
probability that a barrier of $\cM$ appears in 
$[a-\f/3-\bub, a + f/3+\bub]$:
 \[
   \ncln^{*} - \ncln \le \bub^{*}\pub = \T^{\bubxp\txp -1}.
 \]
For sufficiently large $\R_{0}$, we will always have
 $\bub^{**}\pub^{*} < 0.5 \bub^{*}\pub$.
Indeed, this says 
$(\T^{\bubxp\txp - 1})^{\txp} < 0.5 \T^{\bubxp\txp - 1}$, 
which is satisfied if 
$\R \ge 21$.
This implies that if~\eqref{e.strong-ncln-ub} holds for $\cM$ then it also
holds for $\cM^{*}$.
For the inequality~\eqref{e.strong-f-intro},
since the scale-up definition~\eqref{e.new-slb} says 
$\slb^{*} - \slb = \slopeincr \g / \f$, the inequality 
 \[
   \slopeincr \g^{*}/\f^{*} < 0.5(\slb_{0} - \slb^{*})
 \]
will be guaranteed if $\R_{0}$ is large.\footnote{More exactly, if
$\R_{0} \ge 401$.}
 \end{proof}

 \begin{lemma}\label{l.emerg-contrib}
If $\R_{0}$ is sufficiently large then the following holds.
Assume that $\cM=\cM^{k}$ is a mazery.
  \begin{enumerate}[\upshape 1.]

    \item\label{i.emerg-ub}
For any point $x$, the sum, over all $w$, of the
probabilities  that a barrier of the emerging type of rank $\r$ 
and size $w$ starts at $x$ is at most $\p(\r)/2$.

    \item\label{i.emerg-lb}
For the emerging barrier type the fitting emerging holes
satisfy the hole lower bound~\eqref{e.hole-lb}.

  \end{enumerate}
 \end{lemma}
 \begin{proof}
Recall Lemma~\ref{l.emerging}.
Let $n = \Flo{\frac{\g - 5\bub}{3 \bub}}$.
For any point $x$, the expression 
 \[
   2\bub |\Btypes| e^{- (1 - \ncln) n \h(\R^{*})}
 \]
is an upper bound on the sum, over all $w$, of the
probabilities  that an emerging barrier of type $\ang{\g}$ (with rank
$\R^{*}$) starts at $x$. 
We have
 \begin{align*}
               n &> \g / (3 \bub) - 8/3 = \T^{\gxp - \bubxp} / 3 - 8/3,
\\   \h(\R^{*})  &= \aux_{3}\T^{-\txp\hxp},
\\  (1 - \ncln) n \h(\R^{*})
        &> \T^{\gxp - \bubxp - \txp\hxp}\aux_{3}/6 - 1.
 \end{align*}
Due to~\eqref{e.hxp-ub}, this expression grows exponentially in $\R$, and
$e^{-(1 - \ncln) n \h(\R^{*})}$ decreases double exponentially in $\R$.
It follows from~\eqref{e.walltypes-ub} that
its multiplier $2 \bub |\Btypes|$ only grows exponentially in a power of
$\R$.
Hence for large enough $\R_{0}$, the product decreases double
exponentially in $\R$.
So, for sufficiently large $\R_{0}$, claim~\ref{i.emerg-ub} 
follows.\footnote{More exactly, 
we need $\R_{0} \ge 1000$.}

To prove claim~\ref{i.emerg-lb}, let $\ag = \ang{\g}$ be the emerging
barrier type, let $a \le b < c$ and $b-a, c-b  \le 6 \bub^{*}$,
and let $E(a, b, c, \ag)$ be the event
that there is a $d \in [b, c-1]$ such that $(a, d)$ is a
jump, and a hole of type $\ag'$ starts at $d$.
We will be done if we prove
 \begin{equation}\label{e.hole-lb-recall}
   \Prob(E(a,b,c,\ag)) \ge (c-b)^{\hxp} \h(\ag).
 \end{equation}
Let $\cF$ be the event that $a$ is strongly right-clean in $\cM$, 
that $b$ is strongly clean and $b + \g - 4\bub$ is strongly left-clean
in $\cM$ and that no barrier of $\cM$ 
occurs in $[a, b + \bub^{*}]$.
By the definition of emerging holes, $\cF$ implies the event
$E(a, b, c, \ag)$, since $[b, b + \g - 4 \bub]$ will be an emerging hole.
Clearly, 
 \begin{equation}\label{e.emerging-hole}
  \Prob(\cF) \ge 1 - 3 \ncln - 7 \bub^{*} \pub.
 \end{equation}
Lemma~\ref{l.ncln-ub} implies $\ncln < 0.25$, and we have
 \begin{equation*}
  7 \bub^{*}\pub = 7 \T^{\bubxp \txp - 1}.
 \end{equation*}
By~\eqref{e.bubxp-lb}, this is $< 0.1$ if $\R_{0}$ is sufficiently 
large.\footnote{More exactly, we need
$\R_{0} \ge 71$.}
Hence the right-hand side of~\eqref{e.emerging-hole} can be lowerbounded by 
$0.1$.
The required lower bound of~\eqref{e.hole-lb} is
 \[
  \begin{split}
   (c-b)^{\hxp} \h(\ag) &\le (6 \bub^{*})^{\hxp} \h(\R^{*})
   = (6\T^{\txp \bubxp})^{\hxp} \h(\R^{*})
    = \aux_{3} 6^{\hxp} \T^{-\txp\hxp(1 - \bubxp)} < 0.1
  \end{split}
 \]
if $\R_{0}$ is sufficiently large.\footnote{More exactly, we need
$\R_{0} \ge 1803$.}
  \end{proof}

 \begin{lemma}\label{l.compound-contrib}
For sufficiently large $\R_{0}$, the following holds.
Assume that $\cM=\cM^{k}$ is a mazery.
After one operation of forming compound types,
for any rank $\r$ and any point $x$, 
the sum, over all $w$, of the probabilities for the occurrence of a
compound barrier of rank $\r$ and size $w$ at point $x$ is at most 
$\p(\r)\R^{-\aux_{1}/2}$.
 \end{lemma}
 \begin{proof}
Let $\ag_{1},\ag_{2}$ be two types with ranks $\r_{1},\r_{2}$.
Assume without loss of generality that $\r_{1} \le \r_{2}$ and
that $\ag_{1}$ is light: $\r_{1} < \R^{*}=\R^{\txp}$.
With these, according to part~\ref{i.inform.compound} of the scale-up
algorithm, we can form compound barrier types
$\ang{\ag_{1},\ag_{2},i}$, as long as $d_{i} < \f$.
This gives a type of rank $\r_{1}+\r_{2}-i$,
for all $i \le \log_{\lg} \f = \R\fxp$.
The bound~\eqref{e.compound-wall-ub} and the definition of $\p(r)$ 
in~\eqref{e.wall-prob} shows
that the contribution by this term to the sum (over $w$) of
probabilities that a barrier of size $w$ and
rank $\r = \r_{1}+\r_{2}-i$ starts at $x$ is at most
 \begin{equation*}
 (d_{i+1} - d_{i}) \p(\r_{1}) \p(\r_{2}) \le \lg^{i+1}\p(\r_{1})\p(\r_{2}) 
  = \aux_{2}^{2}\lg^{-(\r_{1}+\r_{2}-i-1)}(\r_{1}\r_{2})^{-\aux_{1}}.
 \end{equation*}
Now we have 
$\r_{1}\r_{2} \ge \R\r_{2} \ge (\R/2)(\r_{1}+\r_{2}) \ge \r\R/2$,
hence the above bound reduces to
$\aux_{2}^{2}\lg^{-\r+1}(\r\R/2)^{-\aux_{1}}$.
The total contribution to the sum for rank $\r$ is therefore at most
 \[
  \begin{split}
    &\aux_{2}^{2}\lg^{-\r+1}(\r\R/2)^{-\aux_{1}}
     |\setof{(i,\r_{1}): i \le \R\fxp,\; \r_{1} < \R^{\txp}}|
    \le \aux_{2}^{2} \lg^{-\r+1}(\r\R/2)^{-\aux_{1}} \fxp\R^{\txp+1}
\\  &= \p(\r)\R^{-\aux_{1}/2} 
   c_{2}2^{\aux_{1}}\lg\fxp\R^{-(\aux_{1}/2 - \txp - 1)}
    <  \p(\r)\R^{-\aux_{1}/2},
  \end{split}
 \]
where in the last step we used
 \begin{equation}\label{e.Rlb.compound-contrib}
    \R > (\aux_{2}2^{\aux_{1}}\lg\fxp)^{\frac{1}{\aux_{1}/2 - \txp - 1}},
 \end{equation}
satisfied if $\R_{0} \ge 26$.
 \end{proof}

 \begin{lemma}\label{l.all-wall-ub}
Suppose that each structure $\cM^{i}$ for $i \le k$ is a mazery.
Then inequality~\eqref{e.p()} holds for $\cM^{k+1}$.
\end{lemma}
 \begin{proof}
By Lemma~\ref{l.walltypes-ub}, each rank $\r$ occurs for at most a constant
number $n = \log_{\txp}\frac{2 \txp}{\txp -1}$ values of $k$.
For every such value but possibly the last one, the probability sum
can only be increased as a result of the two
operations of forming compound types.
According to Lemma~\ref{l.compound-contrib}, the increase is upperbounded by
$\p(\r)\R^{-\aux_{1}/2}$.
After these increases, the probability becomes
at most $2 n\p(\r)\R^{-\aux_{1}/2}$.
The last contribution, due to the emerging type,
is at most $\p(\r)/2$ by Lemma~\ref{l.emerg-contrib}; clearly,
if $\R_{0}$ is sufficiently large, the total is still less than 
$\p(\r)$.\footnote{More exactly, we need
$\R_{0} \ge 3$.}
 \end{proof}

 \begin{lemma}\label{l.all-compound-hole-lb}
After choosing $\aux_{1},\aux_{0},\R_{0}$ sufficiently large in this order, 
the following holds.
Assume that $\cM=\cM^{k}$ is a mazery: then every
compound hole type $\bg'$ satisfies the hole lower bound~\eqref{e.hole-lb}.
 \end{lemma}

 \begin{proof}
We will show that compound hole types in $\cM^{*}$
satisfy~\eqref{e.hole-lb} if their component types do (they are either
in $\cM$ or are formed in the process of going from $\cM$ to $\cM^{*}$).
Consider the compound hole type $\bg'$ where
 \[
  \bg = \ang{\ag_{1}, \ag_{2}, i}.
 \]
Let $\r_{j} = \Rank(\ag_{j})$, then $\r = \Rank(\bg) = \r_{1} + \r_{2} - i$.
Let $a \le b < c$ and $b-a, c-b  \le 6 \bub^{*}$.
Following the notation of Lemma~\ref{l.compound-hole-lb},
let $E_{2}(a, b, c, \bg)$ be the event
that there is a $d \in [b, c-1]$ such that $[a,d]$ is a
jump of $\cM^{*}$, and a compound hole of type $\bg'$ starts at $d$.
That lemma assumes $c-a \le 12\bub^{*}$, which holds in our case.
Let us check the condition $(\bub^{*})^{\hxp}\h(\ag_{i}) \le 0.5$.
We have 
 \[
      \h(\ag_{i}) = \aux_{3} \lg^{-\hxp\r_{i}} \le \aux_{3} \T^{-\hxp},
\quad (\bub^{*})^{\hxp} \h(\ag_{i}) \le \aux_{3}\T^{-\hxp(1 - \bubxp\txp)}
 \]
which, due to~\eqref{e.bubxp-lb}, is always smaller than $1/2$ if $\R_{0}$
is sufficiently large.\footnote{More exactly, we need
$\R_{0} \ge 3257$.}
The condition $\pub^{*} \le \pub$ of the lemma is satisfied automatically
by the definitions.
Hence all conditions of the lemma are satisfied.
The conclusion is
 \begin{equation}\label{e.compound-hole-lb-again}
 \begin{split}
   \Prob(E_{2}(a,b,c,\bg)) 
\ge 0.01 (1 - 2\ncln^{*} - 25 \bub^{*}\pub) 
  (c-b)^{\hxp} (e_{i+1} - e_{i})^{\hxp} \h(\ag_{1}) \h(\ag_{2}).
 \end{split}
 \end{equation}
Le us show that for $\aux_{0}$ and then $\R_{0}$ chosen sufficiently large,
this is always larger than $(c-b)^{\hxp}\h(\bg)$.
First we show
 \begin{equation}\label{e.hole-jump-lb}
  e_{i+1}-e_{i} \ge \lg^{i}/17.
 \end{equation}
Indeed, recall the definition of $e_{i}$ in~\eqref{e.hole-jump-types}.
For $i > 17$, we have 
 \[
 e_{i+1}-e_{i} = \flo{\lg^{i}} - \flo{\lg^{i-1}} 
 \ge \lg^{i} - \lg^{i-1} - 1 
= \lg^{i}(1-\lg^{-1}-\lg^{-i}) > 0.1\lg^{i}.
 \]
For $i \le 17$, we have $e_{i+1}-e_{i} \ge 1 \ge \lg^{i}/17$.
This proves~\eqref{e.hole-jump-lb}.
Using~\eqref{e.hole-jump-lb} gives
 \begin{equation}\label{e.all-compound-hole-lb-1}
 \begin{split}
  \h(\ag_{i}) &= \aux_{3}\lg^{-\r_{i}\hxp},
\\ (e_{i+1}-e_{i})^{\hxp}\h(\ag_{1})\h(\ag_{2}) &\ge 17^{-\hxp}
  \aux_{3}^{2}\lg^{-\hxp(\r_{1}+\r_{2} - i)} 
 = 17^{-\hxp}\aux_{3}^{2}\lg^{-\r\hxp}.
 \end{split}
 \end{equation}
Note also that $25\bub^{*}\pub = 25 \T^{\txp\bubxp - 1}< 0.25$ if $\R_{0}$
is large enough.\footnote{More exactly, we need
$\R_{0} \ge 67$.}
Thus, the second factor on the right-hand side of
of~\eqref{e.compound-hole-lb-again} is $\ge 1 - 0.5 - 0.25 = 1/4$.
Substituting into~\eqref{e.compound-hole-lb-again}, we get
the lower bound $\frac{1}{400\cdot 17^{\hxp}}\aux_{3}$ 
for the factors of $(c - b)^{\hxp}\h(\r)$.
This is $\ge 1$ if $\aux_{3}$ is sufficiently large.
More exactly, we need
\begin{equation}\label{e.Rlb.aux3-lb}
  \aux_{3} > 435.5,
 \end{equation}
satisfied by the choice in~\eqref{e.Rlb.aux3}.
 \end{proof}

 \begin{proof}[Proof of Lemma~\protect\ref{l.main}]
The construction of $\cM^{k}$ is complete by the algorithm of
Section~\ref{s.plan}, and the fixing of all parameters in the
present section.

We have to prove that every structure $\cM^{k}$ is a mazery.
The proof is by induction.
We already know that the statement is true for $k=1$: it was 
handled in Example~\ref{x.base}.
Assuming that it is true for all $i \le k$, we prove it for $k+1$.

Condition~\ref{cnd.types}.\ref{i.types.bub} is satisfied by the
argument before Lemma~\ref{l.params-choice-1}.
Condition~\ref{cnd.types}.\ref{i.types.fitting} is satisfied by the
form of the definition of the new types.

Condition~\ref{cnd.distr}.\ref{i.distr.indep} is satisfied as shown in
Lemma~\ref{l.distr.indep}.

Condition~\ref{cnd.distr}.\ref{i.distr.combinat} has been proved in
Lemmas~\ref{l.cover} and~\ref{l.clean}.

In Condition~\ref{cnd.distr}.\ref{i.distr.bounds},
inequality~\eqref{e.p()} has been proved in
Lemma~\ref{l.all-wall-ub}.
Inequality~\eqref{e.0.6} has been proved in
Lemma~\ref{l.0.6}.
Inequality~\eqref{e.ncln-ub} has been proved in Lemma~\ref{l.ncln-ub}.
Inequality~\eqref{e.hole-lb} is proved for emerging walls 
in Lemma~\ref{l.emerg-contrib}, and for compound walls in 
Lemma~\ref{l.all-compound-hole-lb}.

Condition~\ref{cnd.reachable}.\ref{i.reachable.wall} 
is satisfied trivially for the emerging type, (as pointed out in
Lemma~\ref{l.emerging-hole}), and proved for the compound type in
Lemma~\ref{l.compound-hole-fitting}.

Condition~\ref{cnd.reachable}.\ref{i.reachable.hop} is satisfied via 
Lemma~\ref{l.grate} (the Grate Lemma), as discussed after
Lemma~\ref{l.new-hop}. 
There are some conditions on $\f,\g,\bub$ required for this lemma.
Of these, \eqref{e.g-intro} follows from~\eqref{e.Rlb.g-intro.final},
while \eqref{e.f-intro} follows from Lemma~\ref{l.ncln-ub}.

Let us show that the conditions preceding the Main Lemma~\ref{l.main}
hold. 
Condition~\ref{cnd.reachable-hop} is implied by
Condition~\ref{cnd.reachable}.\ref{i.reachable.hop}.
Condition~\ref{cnd.dense} is implied by
Condition~\ref{cnd.distr}.\ref{i.distr.clean}.
Condition~\ref{cnd.Qbd} follows immediately from the definition of
cleanness.

Finally, inequality~\eqref{e.main} follows from~\eqref{e.sum-bub-k+1}.
 \end{proof}

 \section{Proof of Lemma~\protect\ref{l.grate}}\label{s.grate-proof}

Let $b_{\d}(j)$
denote the starting point of wall $W_{\d,j}$ for $j=1,\ldots,n_{\d}$.
Let $b_{\d}(0), b_{\d}(n_{\d}+1)$ denote the beginning and end of interval
$I_{\d}$. 
For convenience, sometimes we will write
 \[
   b(j) = b_{0}(j), \quad c(j) = b_{1}(j).
 \]
Without loss of generality, assume $b_{\d}(0) = 0$.
Let
  \[
  \begin{split}
   L_{\d}(j) &= [b_{\d}(j)+2\g, b_{\d}(j+1)],
\\   P_{i,j} &= (\{b_{0}(i+1)\} \times L_{1}(j)) 
                \cup (L_{0}(i) \times \{b_{1}(j+1)\})
  \end{split}
 \]
for $i=0,\dots,n_{0}$, $j=0,\dots,n_{1}$.
Imagine the rectangle $I_{0} \times I_{1}$
with the direction 0 running horizontally
and the direction 1 vertically, divided into subrectangles by the
vertical lines $x_{0} = b(j)$ ($1 \le j \le n_{0}$), and the
corresponding horizontal lines.
The set $P_{m,n}$ is almost the whole
upper right rim of the $(m,n)$th subrectangle: a segment of size
$2\g$ 
is missing from the left end of the top rim and from the bottom of the
right rim.
For induction purposes, 
we will prove a statement slightly stronger than the lemma.
Let 
 \begin{align*}
    \Gg_{0} &= \bigcup_{m,n}\setof{(x,y)\in P_{m,n} : y-\slb x \ge 8 \g(m+n)},
\\  \Gg_{1} &= \bigcup_{m,n}\setof{(x,y)\in P_{m,n} : x-\slb y \ge 8 \g(m+n)}.
 \end{align*}
Let $H$
be the set of left-clean points $(x_{0}, x_{1})$ (both $x_{0}$ and
$x_{1}$ must be left-clean).
We will show that 
 \[
   (\bigcup_{m,n} P_{m,n}) \cap \Gg_{0} \cap \Gg_{1} \cap H
 \]
is reachable.
Let us first see that this is sufficient.
Note that if $(x,y) \in P_{m,n}$ then $m f / 3  <  x$,
and therefore $m + n < 3 (x + y) / \f$.
Suppose that for $(x,y) \in \bigcup_{m,n} P_{m,n}$ 
with $x \ge y$ we have $(x,y) \not\in \Gg_{0}$.
Then we have 
 \[
   y/x < \slb  + 24 (x+y) \g /(\f x) \le \slb + 48 \g / \f,
 \]
saying that $\minslope(x,y) < \slb^{*}$.

Our claim says that the reachable region is somewhat decreased from
the ``cone'' $\setof{u: \minslope((0,0),u) \ge \slb}$.
Every time the lower side of the reachable region
crosses a vertical line $x = b(i)$ or a horizontal line
$y = c(j)$, it continues in the same direction, but after an upward shift 
by at most $2 \g$.
The upper side gets shifted down similarly.
The two conditions $\Gg_{0}$ and $\Gg_{1}$ are symmetric: the first 
one limits the lower half of the set (where $x \ge y$), the other one the
upper half (where $y \ge x$).
 \begin{figure}


 \includegraphics{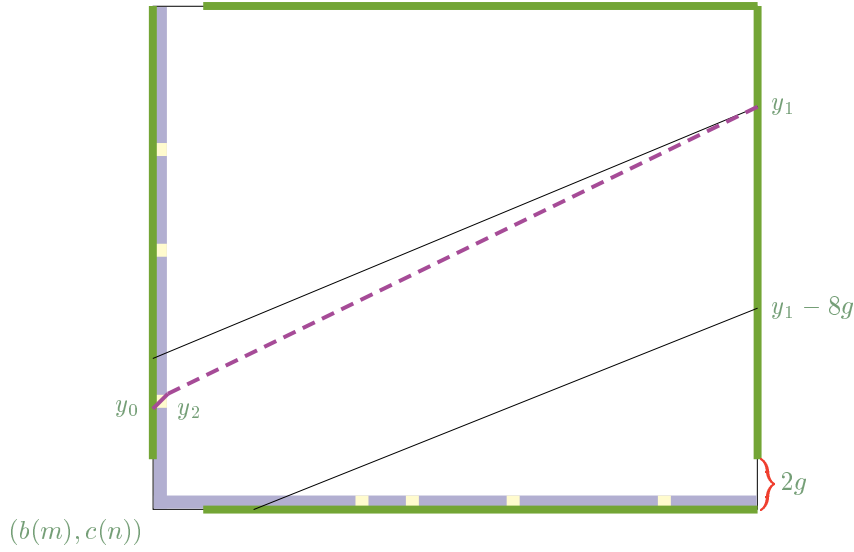}
 \caption{ \label{fig.grate-proof}
To the proof of the Grate Lemma}
 \end{figure}
We will prove the claim by induction on $m+n$; the case $m+n=0$ is
immediate from the reachability
condition~\ref{cnd.reachable}.\ref{i.reachable.hop}.
Consider a point 
 \[
   u_{1} = (x_{1},y_{1})
 \]
in $P_{m,n} \cap \Gg_{0} \cap \Gg_{1} \cap H$.
Without loss of generality, assume $y_{1} \le x_{1}$.
In the interval $b(m) \le x \le b(m+1)$, define the function
 \[
   K(x) = \min(x, y_{1} + (x-x_{1})\slb),
 \]
whose graph is a broken line below the diagonal $x=y$ made up 
of a part (maybe of length 0) of slope 1 followed by a part (maybe of
length 0) of slope $\slb<1$, and ending in point $u_{1}$.
We define the stripe
 \[
 \setof{(x, K(x) - y) : b(m) \le x \le b(m+1),\; 0 \le y \le 8 \g}
 \]
of vertical width $8\g$, below this broken line.
(See Figure~\ref{fig.grate-proof}.)
It intersects the set
$(\{b(m)\} \times L_{1}(n)) \cup (L_{0}(m) \times \{c(n)\})$.
Assume that the intersection with $\{b(m)\} \times L_{1}(n)$ is longer.
The size of this intersection segment is at least $2\g$ (it is not $4\g$
since $L_{1}(n)$ starts only at $c(n)+2\g$, not at $c(n)$.
Its top edge is the point $(b(m), K(b(m)))$.
Its subsegment
 \[
   \{ b(m) \} \times (K(b(m)) + [-1.5\g, -0.5\g])
 \]
contains therefore the starting point
 \[
   u_{0} = (b(m), y_{0})
 \]
of an outer-clean hole fitting the wall $W_{0,m}$.
We have
 \begin{equation}\label{e.y0-bds}
  0.5\g \le K(b(m)) - y_{0} \le  1.5\g.
 \end{equation}
The point $u_{0}$ is left-clean since $b(m)$ is the
start of a wall and $y_{0}$ is the start of an outer-clean hole.
Let 
 \[
  w_{0}, w_{1}
 \]
be the size of the wall $W_{0,m}$ and the size of the hole at $y_{0}$
respectively.
Let us show that $u_{0}$ is reachable.
Since $u_{0}\in P_{m-1,n}$, if 
we show that $u_{0} \in \Gg_{0}\cap \Gg_{1}$ 
then the statement follows from the inductive assumption.
By the definition of our stripe, we have $y_{0} \le b(m)$,
hence $b(m)-\slb y_{0} \ge y_{0} - \slb b(m)$.
Therefore for $u_{0}  \in \Gg_{0} \cap \Gg_{1}$, we only
need to show
 \begin{equation}\label{e.prev-Gg0}
  y_{0} - \slb b(m)  \ge 8 \g(m+n-1).
 \end{equation}
By our assumptions, $u_{1} \in \Gg_{0}$, that is 
$y_{1} - \slb x_{1} \ge 8 \g(m+n)$.
We passed from $u_{0}$ to $u_{1}$ by moving horizontally by the amount
$w_{0}$, moving up
by at most $8 \g$ and then ascending with slope at least $\slb$.
From this, inequality~\eqref{e.prev-Gg0} follows.

Then $u_{0}\leadsto u_{2} = (b(m)+w_{0}, y_{0}+w_{1}) \eqdef (x_{2},y_{2})$.
Also, $x_{2},y_{2}$ are right-clean.
If we show that $\minslope(u_{2},u_{1}) \ge \slb$ then
$u_{2}\leadsto u_{1}$ follows.
This can be done 
using the assumptions $\g > 2\bub \ge 2w_{i}$ and~\eqref{e.y0-bds}:
 \begin{equation*}
  \begin{split}
     y_{1} - K(b(m)) + 0.5\g &\le y_{1} - y_{0} \le y_{1} - K(b(m)) + 1.5\g,
\\       \slb (x_{1} - b(m)) &\le   y_{1} - K(b(m)) \le x_{1} - b(m),
\\  \slb (x_{1} - b(m)) + 0.5\g &\le y_{1} - y_{0} \le x_{1} - b(m) + 1.5\g,
\\             y_{1} - y_{2} &= y_{1} - y_{0} - w_{1},
\\        \slb (x_{1}-x_{2}) &\le \slb (x_{1} - b(m)) + 0.5\g - w_{1} 
\\         &\le y_{1} - y_{2} \le x_{1} - b(m) + 1.5\g - w_{1}
\\    &= x_{1} - x_{2} + 1.5\g - w_{1} + w_{0} \le x_{1} - x_{2} + 2\g.
\\ \sg &\le \frac{y_{1}-y_{2}}{x_{1}-x_{2}} \le 1 + \frac{2\g}{x_{1}-x_{2}} 
  \le 2 \le 1/\slb.
  \end{split}
 \end{equation*}
The case when the larger part of the width of the stripe
intersects the horizontal segment $L_{0}(m) \times \{c(n)\}$, is similar.

\section{Related synchronization problems}\label{s.related}

\begin{figure}
\if\picsize\SlideSize
    \psset{unit=0.5cm}
\else\if\picsize\PaperSize
    \psset{unit=0.8cm}
\fi\fi

\if\palette\DefaultColors
    \definecolor{my-orchid}{cmyk}{0.00, 0.40, 0.02, 0.00}
\else\if\palette\SlideColors
    \definecolor{my-orchid}{cmyk}{0.00, 0.20, 0.01, 0.00}
\else\if\palette\GrayColors
    \definecolor{my-orchid} {gray}{0.8}
\fi\fi\fi

 \[
  \includegraphics{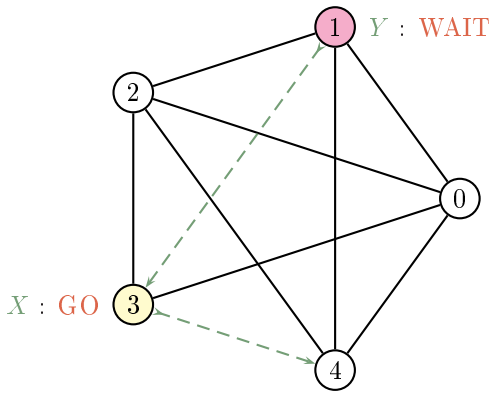}
 \]
 \caption{ \label{fig.rwalk}
The clairvoyant demon problem.
$X,Y$ are ``tokens'' performing independent random walks on the same
graph: here the complete graph $K_{5}$.
A ``demon'' decides every time, whose turn it is.
She is clairvoyant and wants to prevent collision.
}
 \end{figure}

The clairvoyant synchronization problem that has
appeared first has also been introduced by Peter Winkler.
We again have two infinite random sequences $X_{\d}$ for $\d=0,1$
independent from each other.
Now, both of them are random walks on the same graph.
(See Figure~\ref{fig.rwalk}.)
Given delay sequences $t_{\d}$, we say that there is a collision at
$(\d,n,k)$ if $t_{\d}(n) \le t_{1-\d}(k) < t_{\d}(n+1)$ and 
$x_{\d}(k)=x_{1-\d}(n)$.
Here, the delay sequence $t_{\d}$ can be viewed as causing the sequence
$X_{\d}$ to stay in state $x(n)$ between times $t_{\d}(n)$ and
$t_{\d}(n+1)$.
A collision occurs when the two delayed walks enter the
same point of the graph.
This problem, called the clairvoyant demon problem, 
arose originally from a certain leader-election problem in distributed
computing.
Consider the case when the graph is the complete graph of size $m$.
Since we have now just a random walk on a graph, there is no real number
like $\p$ in the compatible sequences problem, that we can decrease 
in order to give a better chance of a solution.
But, $1/m$ serves the same purpose.
Simulations suggest that the walks do not collide if $m \ge 5$, and it
is known that they do collide for $m=3$.
In paper~\cite{GacsWalks01}, we prove that for
sufficiently large $m$, the walks do not collide.
The proof relies substantially on the technique developed here.

\begin{figure}




\[
\includegraphics{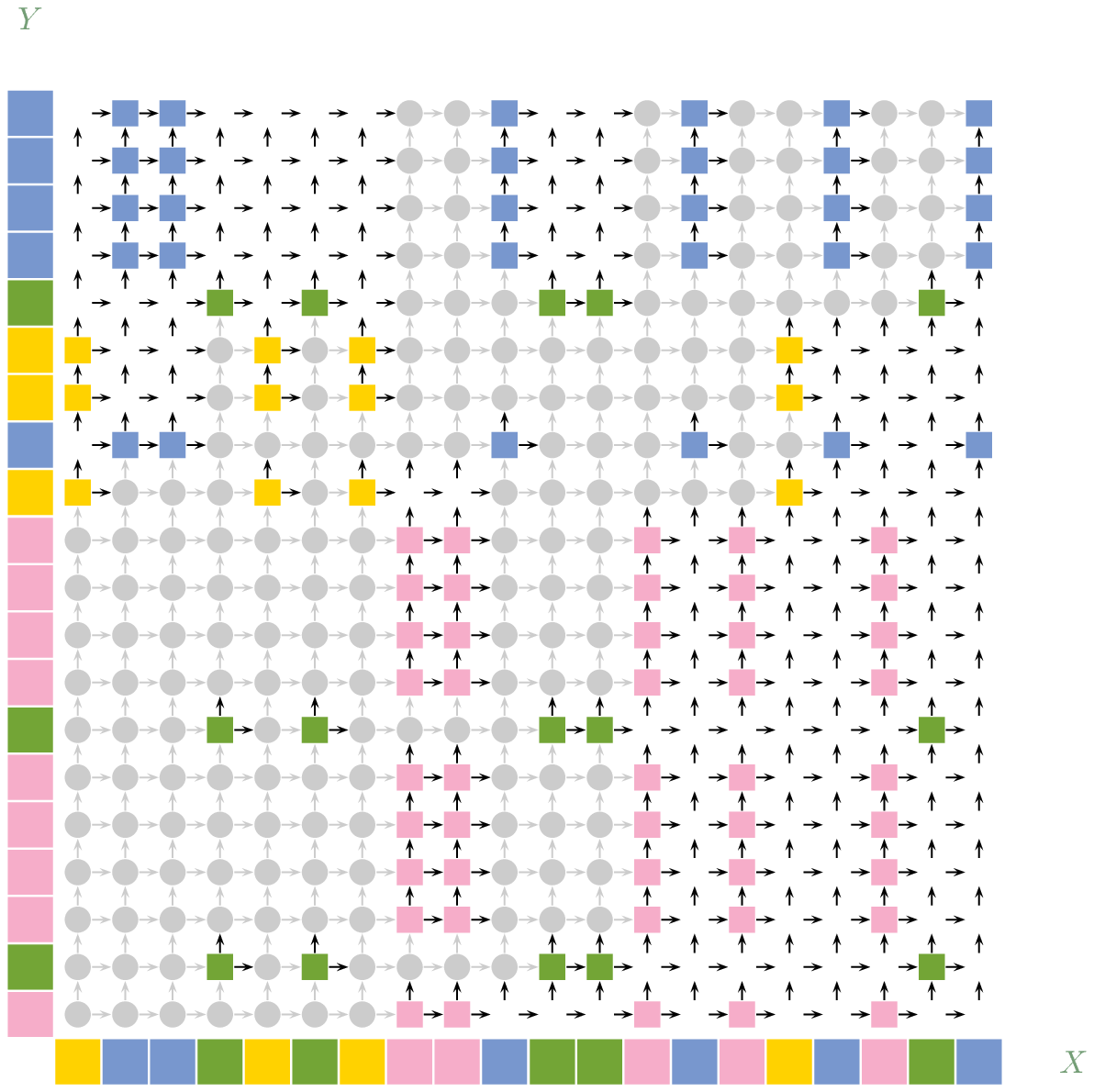}
\]
\[
\mbox{
\includegraphics{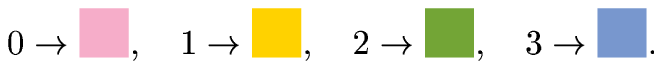}
}
\]
 \caption{
Percolation for the clairvoyant demon problem, for random
walks on the complete graph $K_{4}$.
Round light-grey dots mark the reachable points.}
 \label{fig.cperc}
 \end{figure}

The clairvoyant demon
problem also has a natural translation into a percolation problem,
this time site percolation rather than edge percolation.
(See Figure~\ref{fig.cperc}.)
Consider the lattice $\ZZ_{+}^{2}$, and a graph obtained from it in
which each point is connected to its right and upper neighbor.
For each $i,j$, let us ``color'' the $i$th vertical line 
by the state $X_{0}(i)$, and the $j$th horizontal line by the state
$X_{1}(j)$.
A point $(i,j)$ will be called blocked if $X_{0}(i)=X_{1}(j)$, if its
horizontal and vertical colors coincide.
The question is whether there is, with positive probability, an infinite
directed path (moving only right and up) starting from $(0,0)$ and avoiding
the blocked points.

This problem permits an interesting variation: undirected percolation,
allowing arbitrary paths in the graph, not only directed ones.
This variation has been solved, independently, in~\cite{WinklerCperc99} and
\cite{BalBollobStacCperc99}.
On the other hand, the paper~\cite{GacsClairv99} shows that the
directed problem has a different nature, since if there is
percolation, it has power-law convergence (the undirected percolations have
the usual exponential convergence).

 \section{Conclusions}\label{s.concl}

One of the claims to interest in this dependent percolation
problem is its power-law
behavior over a whole range of parameter values $\p$, not only at a
critical point.
Let us call $b(n)$ the probability that there is a path from $(0,0)$
to distance $n$ but not further.
As we indicated in Subsection~\ref{ss.power-law}, one can prove that
 \[
   b(n) > n^{k_{1}\frac{\log \p}{\p}}
 \]
for some $k_{1}$.
Implicitly, one can see that our paper gives an upper bound 
$n^{-k_{2} / \bubxp} > b(n)$ 
for some $k_{2}$.
Here, $\bubxp$ is clearly not smaller than $\p$, but we do not know 
whether our proof can be refined to make $\bubxp$ approach $\p$.

We could allow $\Prob\setof{X_{\d}(i)=1} = \p_{\d}$ to be different in
the two sequences, say $\p_{0} < \p_{1}$.
This describes the chat
situation when one of the two speakers is more likely to speak than the
other.
It does not seem difficult to
generalize the methods of the present paper to show that
synchronization is possible if $\p_{0}/(1 - \p_{1})$ is small.

\subsection*{Acknowledgement}

The author is grateful to Peter Winkler, Marton Balazs and John Tromp 
for valuable comments, 
and particularly the anonymous referee for his very careful reading that
caught many errors.

\end{document}